\documentclass[final,onefignum,onetabnum]{siamart171218}

\usepackage{amsfonts}
\usepackage[english]{babel}
\usepackage[utf8]{inputenc}
\usepackage{amsmath}
\usepackage{graphicx}
\usepackage{amssymb}
\usepackage{mathrsfs}
\usepackage[colorinlistoftodos]{todonotes}
\usepackage{mathtools}
\usepackage{bm}
\usepackage{adjustbox}   
\usepackage{subfigure}
\usepackage[mmddyyyy]{datetime}

\usepackage{tabularx,colortbl}
\usepackage{colortbl}
\usepackage{graphicx}

\colorlet{mylinkcolor}{blue}
\colorlet{mycitecolor}{orange}
\colorlet{myurlcolor}{orange}

\newcommand{\vcr}[1]{\bm{#1}}
\newcommand{\mat}[1]{\bm{#1}}
\newcommand{\tsr}[1]{\pmb{\mathcal{#1}}}

\newcommand{\fnrm}[1]{{\| #1 \|}_F}

\newcommand{\inti}[2]{\{{#1},\ldots, {#2}\}}

\newcommand{\R}{\mathbb{R}}

\usepackage[normalem]{ulem} 

\usepackage{cite}
\usepackage{amsmath,amssymb,amsfonts}

\usepackage{algpseudocode}
\usepackage{algorithm}

\usepackage{graphicx}
\usepackage{textcomp}
\usepackage{xcolor}

\usepackage{listings}
\definecolor{mygreen}{rgb}{0,0.2,0}
\definecolor{mygray}{rgb}{0.5,0.5,0.5}
\definecolor{mymauve}{rgb}{0.58,0,0.82}
\definecolor{mypurple}{rgb}{0.38,0,0.32}
\definecolor{myblue}{rgb}{0.1,0,0.32}
\newcommand{\costyle}{\footnotesize\ttfamily\bfseries}
\newcommand{\kwstyle}{\costyle\textcolor{myblue}}
\def\BibTeX{{\rm B\kern-.05em{\sc i\kern-.025em b}\kern-.08em
    T\kern-.1667em\lower.7ex\hbox{E}\kern-.125emX}}
\begin{document}





\title{Comparison of Accuracy and Scalability of Gauss-Newton and Alternating Least Squares for CP Decomposition
}

\author{Navjot Singh\thanks{Department of Mathematics, University of Illinois at Urbana-Champaign, Urbana, IL, 61801 (\email{navjot2@illinois.edu}).} \and Linjian Ma\thanks{Department of Computer Science, University of Illinois at Urbana-Champaign, Urbana, IL, 61801 (\email{lma16@illinois.edu}, \email{hyang87@illinois.edu}, \email{solomon2@illinois.edu}).}
\and Hongru Yang\footnotemark[3] \and Edgar Solomonik\footnotemark[3]
}

\maketitle

\begin{abstract}
Alternating least squares is the most widely used algorithm for CP tensor decomposition. However, alternating least squares may exhibit slow or no convergence, especially when high accuracy is required. An alternative approach is to regard CP decomposition as a nonlinear least squares problem and employ Newton-like methods. Direct solution of linear systems involving an approximated Hessian is generally expensive. However, recent advancements have shown that use of an implicit representation of the linear system makes these methods competitive with alternating least squares. We provide the first parallel implementation of a Gauss-Newton method for CP decomposition, which iteratively solves linear least squares problems at each Gauss-Newton step. In particular, we leverage a formulation that employs tensor contractions for implicit matrix-vector products within the conjugate gradient method. The use of tensor contractions enables us to employ the Cyclops library for distributed-memory tensor computations to parallelize the Gauss-Newton approach with a high-level Python implementation. 
In addition, we propose a regularization scheme for Gauss-Newton method to improve convergence properties without any additional cost.
We study the convergence of variants of the Gauss-Newton method relative to ALS for finding exact CP decompositions as well as approximate decompositions of real-world tensors. We evaluate the performance of sequential and parallel versions of both approaches, and study the parallel scalability on the Stampede2 supercomputer. 

\end{abstract}



\ifpdf
\hypersetup{ pdftitle={gauss newton cpd} }
\fi

\headers{Gauss-Newton for CP Decomposition}{Navjot Singh, Linjian Ma, Hongru Yang and Edgar Solomonik}

\section{Introduction}
\label{sec:intro}
The CP (canonical polyadic or CANDECOMC/PARAFAC) tensor decomposition
is widely used for data analytics in different scientific fields~\cite{murphy2013fluorescence,sidiropoulos2017tensor,maruhashi2011multiaspectforensics,cong2015tensor,hitchcock1927expression }, machine learning applications~\cite{anandkumar2014tensor,kolda2009tensor,bailey2017word}, and quantum chemistry~\cite{thomas2017intertwined}. 
CP decomposition of an input tensor can be computed via different optimization techniques, such as variants of gradient descent~\cite{acar2011scalable,paatero1997weighted}, deflations~\cite{anandkumar2014guaranteed,anandkumar2014tensor}, and alternating least squares~\cite{kolda2009tensor}.

Nowadays, the alternating least squares (ALS) method, which solves quadratic optimization subproblems for each factor matrix in an alternating manner, is most commonly used and has become a target for parallelization~\cite{karlsson2016parallel,hayashi2017shared}, performance optimization~\cite{ma2018accelerating,schatz2014exploiting}, and acceleration by randomization~\cite{battaglino2017practical}.
A major advantage of ALS is its guaranteed monotonic decrease of the residual.
However, there are many cases where ALS shows slow or no convergence when solution with high resolution is required, which is also called the `swamp' phenomenon~\cite{mitchell1994slowly}.
Swamps deteriorate both the running time and the convergence behavior of the ALS method.
Consequently,
researchers have been looking at different alternatives to ALS, including various regularization techniques~\cite{navasca2008swamp,li2013some}, line search~\cite{rajih2008enhanced,nion2008enhanced,mitchell2018nesterov} and gradient based methods~\cite{acar2011scalable,paatero1997weighted,phan2013low,sorber2013optimization,tichavsky2013further,tomasi2006comparison}.

Of the variants of gradient based methods, one promising approach is to perform the CP decomposition by solving a nonlinear least squares problem using the Newton or Gauss-Newton methods~\cite{paatero1997weighted,tomasi2006comparison,tomasi2005parafac}.
These approaches offer superlinear convergence and are better at avoiding the swamps inhibiting performance of ALS.
Naive solution of linear equations arising in these methods is expensive to get.
For rank-$R$ decomposition of an order-$N$ tensor with all the dimension sizes equal to $s$, standard algorithms either perform Cholesky on the normal equations~\cite{paatero1997weighted} or QR on the Jacobian matrix~\cite{tomasi2006comparison}, yielding a complexity of $O(N^3s^3R^3)$.
However, the matrices involved in this linear system are sparse and have much implicit structure.
A recent advancement has shown that the cost of inverting the Hessian can be reduced to $O(N^3R^6)$~\cite{phan2013low}. A successive study showed that the cost can be further reduced to $O(NR^6)$, albeit the approach can suffer from numerical instability \cite{tichavsky2013further}. 
 
Another approach for performing Gauss-Newton with low cost is to leverage an implicit conjugate gradient (CG) method~\cite{sorber2013optimization}.
The structure of the approximated Hessian can be leveraged to perform fast matrix-vector multiplications for CG iterations (with a cost of $O(N^2sR^2)$ per iteration), an approach that can also be augmented with preconditioning to accelerate CG convergence rate~\cite{sorber2013optimization}.
In comparison to the aforementioned direct methods, this iterative approach is substantially more scalable with respect to the CP rank $R$.
This advantage is critical in many applications of CP decomposition, as in many cases $R\geq s$ is needed (in general CP rank can be as high as $s^{N-1}$ for an order $N$ tensor).
Moreover, for the CP decomposition with rank $R< s$, Tucker decomposition (or simply HoSVD)~\cite{tucker1966some} can be used to effectively compress the input tensor from dimensions of size $s$ to $R$, and then CP decomposition can be performed.

In this paper, we investigate the behavior of Gauss-Newton optimization with preconditioned CG on CP decomposition in high rank scenarios (with $R\geq s$ or more generally when the rank is at least the smallest dimension size of the input tensor). 
We consider various approaches to regularization for Gauss-Newton with implicit CG and ALS.
To understand their efficacy, we quantify their ability to converge to exact CP decompositions of synthetic tensors of various CP ranks, as well as to approximate tensors arising in applications in quantum chemistry.
With the best regularization strategy, we find that Gauss-Newton is able to consistently find exact CP decompositions for problems where ALS generally does not converge.
Further, the Gauss-Newton method obtains lower residuals in approximation. We present these results in Section~\ref{sec:exp}.

Our main contribution is the parallel implementation of Gauss-Newton with implicit CG, via a tensor-contraction-based formulation of the method.
We develop a distributed-memory implementation of the method using the Cyclops library for parallel tensor algebra.
Our implementation supports both NumPy and Cyclops as backends, enabling both sequential and parallel experimental studies. We detail our implementations in Section~\ref{sec:implement}.
We evaluate the strong and weak scalability of the method on the Stampede2 supercomputer, and compare its performance to ALS for a variety of test problems.
Our results demonstrate that the Gauss-Newton method can converge faster both in sequential and parallel settings. These results are presented in Section~\ref{sec:exp}.

This paper makes the following contributions:

\begin{itemize}

    \item 
We cast the large matrix-vector multiplication into several tensor contractions so that an existing library on parallel tensor contractions can be utilized. Our analysis achieves the same computational cost as previous work~\cite{sorber2013optimization}.
    \item We propose and evaluate a new regularization strategy, and demonstrate that it is well-suited for CP decomposition with Gauss-Newton.
    \item We provide the first parallel implementation of Gauss-Newton for CP decomposition.
    \item We demonstrate that an implementation of parallel Gauss-Newton with preconditioned CG can both converge faster and achieve higher convergence probability for CP decompositions of both synthetic and application-based tensors with high CP rank.
\end{itemize}

\section{Background}
\label{sec:bg}
We introduce the notation and definitions used in the forthcoming sections here along with a brief introduction to the alternating least squares algorithm.
We suggest~\cite{vannieuwenhoven2015computing,kaya2016parallel,ballard2018parallel, kolda2009tensor,ma2018accelerating} for a detailed review of the algorithm and it's HPC formulation. 

\subsection{Notation and Definitions}

We use tensor algebra notation in both element-wise form and specialized form  for tensor operations~\cite{kolda2009tensor}.
For vectors, bold lowercase Roman letters are used, e.g., $\vcr{x}$. For matrices, bold uppercase Roman letters are used, e.g., $\mat{X}$. For tensors, bold calligraphic fonts are used, e.g., $\tsr{X}$. An order $N$ tensor corresponds to an $N$-dimensional array with dimensions $s_1\times \cdots \times s_N$. 
Elements of vectors, matrices, and tensors are denoted in subscript, e.g., $x_i$ for a vector $\vcr{x}$, $x_{ij}$ for a matrix $\mat{X}$, and $x_{ijkl}$ for an order 4 tensor $\tsr{X}$. 
The $i$th column of a matrix $\mat{X}$ is denoted by $\vcr{x}_i$. 
The mode-$n$ matrix product of a tensor $\tsr{X} \in \mathbb{R}^{s_1 \times \cdots \times s_N}$ with a matrix $ \mat{A}\in \mathbb{R}^{J\times s_n}$ is denoted by $\tsr{X}\times_n \mat{A}$, with the result having dimensions $s_1\times\cdots\times s_{n-1}\times J\times s_{n+1}\times\cdots\times s_N$. 
Matricization is the process of reshaping a tensor into a matrix. Given a tensor $\tsr{X}$ the mode-$n$ matricized version is denoted by $\mat{X}_{(n)}\in \mathbb{R}^{s_n\times K}$ where $K=\prod_{m=1,m\neq n}^N s_m$. 
We use parenthesized superscripts as labels for different tensors and matrices, e.g., $\mat{A}^{(1)}$ and $\mat{A}^{(2)}$ are different matrices.

The Hadamard product of two matrices $\mat{U}, \mat{V} \in \mathbb{R}^{I\times J}$ resulting in matrix $\mat{W} \in \mathbb{R}^{I\times J}$ is denoted by $\mat{W} = \mat{U} \ast \mat{V}$, where $w_{ij}= u_{ij}v_{ij}$. 
The inner product of matrices $\mat{U}, \mat{V}$ is denoted by $\langle \mat{U}, \mat{V} \rangle = \sum_{i,j}u_{ij}v_{ij}$.
The outer product of K vectors $\vcr{u}^{(1)}, \ldots , \vcr{u}^{(K)}$ of corresponding sizes $s_1, \ldots , s_K$ is denoted by $\tsr{X} = \vcr{u}^{(1)} \circ \cdots \circ \vcr{u}^{(K)}$ where $\tsr{X} \in \mathbb{R}^{s_1 \times \cdots \times s_K}$ is an order $K$ tensor. 
%
For matrices $\mat{A}\in \mathbb{R}^{I\times K} = \begin{bmatrix} \vcr{a}_1, \ldots, \vcr{a}_K \end{bmatrix}$ and $\mat{B}\in \mathbb{R}^{J\times K}= \begin{bmatrix} \vcr{b}_1, \ldots, \vcr{b}_K \end{bmatrix}$, their Khatri-Rao product resulting in a matrix of size $(IJ)\times K$ defined by
   $
   \mat{A}\odot \mat{B} = [\vcr{a}_1\otimes \vcr{b}_1,\ldots, \vcr{a}_K\otimes \vcr{b}_K]
   $, where $\vcr{a}\otimes \vcr{b}$ denotes the Kronecker product of the two vectors.

\subsection{CP Decomposition with Alternating Least Squares}

The CP tensor decomposition~\cite{hitchcock1927expression,harshman1970foundations} for an input tensor $\tsr{X}\in \R^{s_1\times\dots\times s_N}$
is denoted by
   \[
    \tsr{X} \approx [\![ \mat{A}^{(1)}, \cdots , \mat{A}^{(N)} ]\!], \quad \text{where} \quad  
    \mat{A}^{(i)} = [ \mat{a}_1^{(i)}, \cdots , \mat{a}_r^{(i)} ],
    \]
and serves to approximate a tensor by a sum of $R$ tensor products of vectors,
  \[
  \tsr{X} \approx \sum_{r=1}^{R} \mat{a}_r^{(1)}\circ \cdots \circ \mat{a}_r^{(N)}.
  \]
The CP-ALS method alternates among quadratic optimization problems for each of the factor matrices $\mat{A}^{(n)}$, resulting in linear least squares problems for each row,
\[
 \mat{A}^{(n)}_{\text{new}}\mat{P}^{(n)}{}^T \cong \mat{X}_{(n)},
\]
where the matrix $\mat{P}^{(n)}\in \mathbb{R}^{I_n \times R}$, where $I_n = s_1\times \cdots \times s_{n-1}\times s_{n+1}\times \cdots \times s_N$ is formed by Khatri-Rao products of the other factor matrices,
\begin{align}
    \mat{P}^{(n)}=\mat{A}^{(1)} \odot \cdots \odot  \mat{A}^{(n-1)}  \odot  \mat{A}^{(n+1)} \odot \cdots \odot \mat{A}^{(N)}. \label{eq:P}
\end{align}
These linear least squares problems are often solved via the normal equations~\cite{kolda2009tensor},
     \[
     \mat{A}^{(n)}_{\text{new}}\boldsymbol{\Gamma}^{(n)}\leftarrow \mat{X}_{(n)}\mat{P}^{(n)},
     \]  
where $\mat{\Gamma}\in\R^{R\times R}$ can be computed via
\begin{equation}
\boldsymbol{\Gamma}^{(n)}=\mat{S}^{(1)}\ast\cdots\ast \mat{S}^{(n-1)}  \ast \mat{S}^{(n+1)}\ast\cdots\ast \mat{S}^{(N)},
\label{eq:gamma}
\end{equation}
with each
$\mat{S}^{(i)} = \mat{A}^{(i)T}\mat{A}^{(i)}.$
These equations also give the $n$th component of the optimality conditions for the unconstrained minimization of the nonlinear objective function,
   \begin{equation}
f(\mat{A}^{(1)}, \ldots , \mat{A}^{(N)}) :=
 \frac{1}{2}\vert \vert\pmb{\mathcal{X}}-[\![ \mat{A}^{(1)}, \cdots , \mat{A}^{(N)} ]\!]\vert \vert_F^2.
 \label{eq:obj}
   \end{equation}
The \textit{Matricized Tensor Times Khatri-Rao Product} or \mbox{MTTKRP} computation $\mat{M}^{(n)}=\mat{X}_{(n)}\mat{P}^{(n)}$ is the main computational bottleneck of CP-ALS\cite{ballard2018communication}. A work efficient way to compute MTTKRP is to contract the factor matrices with the tensor successively. The bottleneck for this implementation is the contraction between the tensor and the first-contracted matrix. Algebraically, this contraction can be written as the tensor times matrix product, $\tsr{X}\times_i \mat{A}^{(i)T}$.
For a rank-R CP decomposition, this computation has the cost of $2s^NR$ if $s_n=s$ for all $n\in\inti{1}{N}$.

The dimension-tree algorithm for ALS~\cite{phan2013fast,vannieuwenhoven2015computing,kaya2016parallel,ballard2018parallel,kaya2017high,kaya2019computing} uses a fixed amortization scheme to update \mbox{MTTKRP} in each ALS sweep. 
This scheme only needs to perform two bottleneck contraction calculations for each ALS sweep, decreasing the leading order cost of a sweep from $O(Ns^{N}R)$ to $ O(s^{N}R). $


\section{Gauss-Newton for CP Decomposition}
\label{sec:gnm}
The Gauss-Newton (GN) method is a modification of the Newton's method to solve nonlinear least squares problem for a quadratic objective function defined as 
\[
\phi(\vcr{x}) = \frac{1}{2}\|\vcr{y}- \vcr{f}(\vcr{x})\|^2,
\]
where $\vcr{y}$ is the given vector of points with respect to which we solve the least squares problem, $\vcr{x}$ is the solution vector required and $\vcr{f}$ is the nonlinear function of $\vcr{x}$ given in the problem. The gradient and the Hessian matrix of $\phi(\vcr{x})$ can be expressed as
\begin{align*}
\nabla \phi (\vcr{x}) = & \mat{J}^T_r(\vcr{x})\vcr{r}(\vcr{x}),  \\
\mat{H}_{\phi}(\vcr{x}) = \mat{J}_r^T(\vcr{x})\mat{J}_r(\vcr{x})& + \sum_{i}r_i(\vcr{x})\mat{H}_{r_i}(\vcr{x}),
\end{align*}
where $\vcr{r}(\vcr{x})$ is the residual function defined as $\vcr{r}(\vcr{x}) = \vcr{y}- \vcr{f}(\vcr{x}) $, $\mat{J}_r(\vcr{x})$ is the Jacobian matrix of the residual function with respect to $\vcr{x}$, and $\mat{H}_{r_i}(\vcr{x})$ is the Hessian matrix of a component of the residual function $r_i$ with respect to $\vcr{x}$.
 
The Gauss-Newton method leverages the fact that $\mat{H}_{r_i}(\vcr{x})$ is small in norm when the residual is small, to approximate the Hessian as $\mat{H}_{\phi}(\vcr{x}) \approx \mat{J}_r^T(\vcr{x}) \mat{J}_r(\vcr{x})$.
Consequently, the Gauss-Newton iteration aims to perform the update,
\[\vcr{x}^{(k+1)} = \vcr{x}^{(k)}- (\mat{J}_r^T(\vcr{x}^{(k)})\mat{J}_r(\vcr{x}^{(k)}))^{-1} \mat{J}_r^T(\vcr{x}^{(k)})\vcr{r}(\vcr{x}^{(k)}),\]
where $\vcr{x}^{(k)}$ represents the $\vcr{x}$ at $k$th iteration.
This linear system corresponds to the normal equations for the linear least squares problem,
\[
\mat{J_r}(\vcr{x}^{(k)})(\vcr{x}^{(k+1)} -\vcr{x}^{(k)}) \cong  -\vcr{r}(\vcr{x}^{(k)}).
\]
Approximation with CP decomposition \eqref{eq:obj} is a nonlinear least squares problem where the points are tensor entries and the unknowns are factor matrix entries.

We define the Jacobian tensor as 
\[
\tsr{J} = [\tsr{J}^{(1)}, \ldots,  \tsr{J}^{(N)}]
\]
 for the N-dimensional CP decomposition, where  $\tsr{J}^{(n)}\in\mathbb{R}^{s_1\times\cdots\times s_N\times s_n \times R}$ is the Jacobian tensor for the residual tensor with respect to $\mat{A}^{(n)}$, and is expressed element-wise as
\begin{equation}
j^{(n)}_{i_1 \dots i_N k r} = \displaystyle{\bigg(-\prod_{m=1,m \neq n}^N a^{(m)}_{i_mr} \bigg) \delta_{i_nk} } .
\label{eq:jacobian}
\end{equation}
Another way to derive the Jacobian matrices is by unfolding the factor matrices and the residual function as suggested in \cite{acar2011scalable}.
Factorization of the Hessian to solve a linear system in Gauss-Newton has a cost of $O(N^3s^3R^3)$.
More advanced approaches to solving Hessian can achieve a cost of $O(NR^6)$~\cite{tichavsky2013further}, but this reduction is not substantial when CP rank is high, i.e., $R\geq s$.

Alternatively, conjugate gradient (CG) with implicit matrix products can be used to solve the linear least squares problems in this Gauss-Newton method~\cite{sorber2013optimization}.
Instead of performing a factorization or inversion of the approximate Hessian matrix, this approach only needs to perform matrix vector products $\mat{J}^T\mat{J} \vcr{v}$ at each iteration (henceforth we drop the subscript $r$ from $\mat{J}_r$ and simply refer to $\mat{J}$ for the matrix form of the Jacobian and $\tsr{J}$ for its tensor form).
We derive the matrix vector product in terms of tensor contractions in the following section.

\subsection{Gauss-Newton with Implicit Conjugate Gradient}
\label{subsec:gn}
With the Jacobian tensors defined in \eqref{eq:jacobian}, the matrix-matrix product $\mat{H} = \mat{J}^T\mat{J}$ can be expressed as an operator with the following form, 
\[
h^{(n,p)}_{krlz} = \displaystyle{ \sum_{i_1 \dots i_N}j^{(n)}_{i_1 \dots i_N kr} j^{(p)}_{i_1 ... i_N lz}},
\]
which can be simplified to
\begin{align}
h^{(n,p)}_{krlz} &= \begin{cases}
      \displaystyle{\delta_{kl} \Gamma_{rz}^{(n,n)}}, & \text{if}\ n=p \\
      \displaystyle{a^{(n)}_{kz}a^{(p)}_{lr}\Gamma_{rz}^{(n,p)}}, & \text{otherwise}
    \end{cases},
\label{eq:J} \\
%
\text{where} \quad \Gamma_{rz}^{(n,p)} &=
\prod_{m=1, m\neq n,p}^N \Bigg(\sum_{i_m}a^{(m)}_{i_mr}a^{(m)}_{i_mz} \Bigg).
\label{eq:gamma2}
\end{align}
Note that $\mat{\Gamma}^{(n,n)}=\mat{\Gamma}^{(n)}$ as defined in \eqref{eq:gamma}.
The matrix-vector product $\mat{H}\vcr{w}$ can be written as
\[
\mat{H}\vcr{w} = \displaystyle{\sum_{n=1}^N \sum_{p=1}^N \sum_{l=1}^{s_p}\sum_{z=1}^R h^{(n,p)}_{krlz}w_{lz}^{(p)}}.
\]
The contractions in the innermost summation have the form,

\begin{equation}
\displaystyle{ \sum_{l,z} h^{(n,p)}_{krlz}w_{lz}^{(p)} }  = \begin{cases}
      \displaystyle{\sum_{z}\,  \Gamma_{rz}^{(n,n)}w_{kz}^{(n)}  },   & \text{if}\ n=p,  \\
      
      \displaystyle{\sum_{l,z}\,  a^{(n)}_{kz}a^{(p)}_{lr}\Gamma_{rz}^{(n,p)}w_{lz}^{(p)}
      }, & \text{otherwise.}
    \end{cases}
    \label{eq:matvecprod}
\end{equation}
Computation of $\displaystyle{\sum_{n=1}^N \sum_{p=1}^N \sum_{l,z} h^{(n,p)}_{krlz}w_{lz}^{(p)}}$ requires $N^2$ contractions of  the form $\displaystyle{\sum_{l,z} h^{(n,p)}_{krlz}w_{lz}^{(p)}}$ for a total cost of 
$O(N^2sR^2)$ when each mode of the input tensor has size $s$ and is $O\Big(N \displaystyle{(\sum_{m=1}^N s_m )R^2}\Big)$ in the general case. Our Gauss-Newton algorithm is summarized in Algorithm~\ref{alg:nls}.

\begin{algorithm}[t]
    \caption{\textbf{CP-GN}: Gauss-Newton with preconditioned implicit CG for CP decomposition}
\label{alg:nls}
\begin{algorithmic}[1]
\small
\State{\textbf{Input: }Tensor $\tsr{X}\in\mathbb{R}^{s_1\times\cdots\times s_N}$, 
stopping criteria $\varepsilon$, CG stopping criteria $\varepsilon_{cg}$, rank $R$}
\State{Initialize $\{\mat{A}^{(1)}, \ldots , \mat{A}^{(N)}\}$ so each $\mat{A}^{(n)}\in\mathbb{R}^{s_n\times R}$ is random

}
\While{$\sum_{i=1}^{N}{\fnrm{\mat{G}^{(i)}}}>\varepsilon$}

\State{Calculate $\mat{M}^{(n)}=\mat{X}_{(n)}\mat{P}^{(n)}$ for $n \in \inti{1}{N} $}

\Statex{\Comment{Using dimension tree with $\mat{P}^{(n)}$ is defined as in \eqref{eq:P}}}

\For{\texttt{$n \in \inti{1}{N} $}}


\State{Calculate $\boldsymbol{\Gamma}^{(n,p)}$ for $p \in \inti{1}{N} $ based on \eqref{eq:gamma2}
}

\State{$\mat{G}^{(n)} \leftarrow \mat{A}^{(n)}\mat{\Gamma}^{(n,n)} - \mat{M}^{(n)}$
}
\EndFor

\State{Define $\lambda$ based on varying scheme described in Section~\ref{subsec:regu}}

\State{$\begin{aligned}
\{\mat{V}^{(1)},\ldots,\mat{V}^{(N)}\} \leftarrow \text{CP-CG} (&\tsr{X}, 
                                                                 \{\mat{G}^{(1)},\ldots,\mat G^{(N)}\}, \\
                                                                &\{\mat{A}^{(1)},\ldots,\mat A^{(N)}\}, \\
                                                                &\{\mat{\Gamma}^{(n,p)} : n,p\in\{1,\ldots,N\}\}, \\
                                                                &\varepsilon_{cg}, \lambda)
\end{aligned}$}

\For{\texttt{$n \in \inti{1}{N} $}}
\State{$\mat{A}^{(n)} \leftarrow \mat{A}^{(n)}+ \mat{V}^{(n)}$
}
\EndFor

\EndWhile

\State{\Return factor matrices $\{\mat{A}^{(1)}, \ldots , \mat{A}^{(N)}\}$ with $\mat{A}^{(n)}\in \mathbb{R}^{s_n\times R}$}
\end{algorithmic}

\end{algorithm}

\subsection{Regularization for Gauss-Newton}
\label{subsec:regu}

Since the approximated Hessian is inherently rank-deficient \cite{tomasi2006comparison}, we incorporate Tikhonov regularization when solving the linear system, $\mat{J}^T\mat{J} + \lambda\mat{I}$, at each iteration, which corresponds to the Levenberg-Marquardt algorithm~\cite{more1978levenberg}. 
The convergence behavior of the Gauss-Newton method for CP decomposition as well as the CG method used within each Gauss-Newton iteration is sensitive to the choice of regularization parameter.

A common approach to resolve the scaling indeterminacy for the linear least squares problem
is to use $\mat{J}^T\mat{J} + \lambda \text{diag}(\mat{J}^T\mat{J})$, 
however, this may not be the best way to regularize as mentioned in~\cite{more1978levenberg} and we observe this case with a constant $\lambda$ parameter. 
There are several other approaches for choosing the damping parameter and the diagonal matrix at each iteration to ensure local convergence of the algorithm~\cite{more1978levenberg}, but they require additional objective function or gradient calculations, which are costly
in the context of CP decomposition, due to the high computational and communication cost associated with each iteration.

We provide a new heuristic for choosing the damping parameter by varying the regularization at each step.
Variable regularization has been used in the past for the Gauss-Newton method, by increasing or decreasing the parameter depending on the value of the objective function at the next iteration~\cite{more1978levenberg}.
We find that for CP decomposition, variation of the regularization parameter is useful for getting out of swamps, and adjusting it eagerly helps avoiding the need for expensive recomputation of the objective function.

In particular, 
we define an upper threshold and a lower threshold, and initialize $\lambda$ near the upper threshold.
This large value ensures that we take steps towards the negative gradient direction, and enables CG to converge quickly.
 Next, we choose a constant hyper parameter $\mu >1$ and update the $\lambda$ at each iteration with $\lambda = \lambda/\mu$. This update is continued until $\lambda$ reaches the lower threshold, and then it is increased by the update $\lambda = \lambda \mu$ until it reaches the upper threshold value and then decreased again. The lower threshold ensures that the conditioning of $\mat{J}^T\mat{J}$ does not affect the CG updates. 

We show in Section~\ref{subsec:convlikelihood} that this type of varying regularization can significantly improve the convergence probability of Gauss-Newton method
relative to a fixed regularization parameter when an exact CP decomposition exists.
We find that this strategy is robust in speed, accuracy and probability of convergence to global minima across many experiments.

\begin{algorithm}[t]
    \caption{\textbf{CP-CG}: Preconditioned implicit CG for CP decomposition}
\label{alg:cg}
\begin{algorithmic}[1]
\small
\State{\textbf{Input: }Tensor $\tsr{X}\in\mathbb{R}^{s_1\times\cdots\times s_N}$,  gradient set $\{\mat{G}^{(1)},\ldots,\mat G^{(N)}\}$, factor matrix set $\{\mat{A}^{(1)},\ldots,\mat A^{(N)}\}$, set of $R\times R$ matrices $\{\mat{\Gamma}^{(n,p)} : n,p\in\{1,\ldots,N\}\}$, stopping criteria $\varepsilon_{cg}$, regularization term $\lambda$}
\For{\texttt{$n \in \inti{1}{N} $}}

\State{
$\mat{P}_{\text{inv}}^{(n)}\leftarrow (\mat{\Gamma}^{(n,n)} + \lambda \mat{I})^{-1}$
}

\State{Initialize $\mat{V}^{(n)}$ \text{to zeros}}

\State{$\mat{R}^{(n)} \leftarrow - \mat{G}^{(n)}$
}

\State{$\mat{Z}^{(n)}\leftarrow \mat{R}^{(n)}\mat{P}_{\text{inv}}^{(n)}$ 
}

\State{$\mat{W}^{(n)}\leftarrow \mat{Z}^{(n)}$}

\EndFor

\While{$\sum_{i=1}^{N}{\fnrm{\mat{R}^{(i)}}}>\varepsilon_{cg} \sum_{i=1}^{N}{\fnrm{\mat{G}^{(i)}}}$}
\For{\texttt{$n \in \inti{1}{N} $}}

\Statex{\Comment{Using implicit matrix vector product as in \eqref{eq:matvecprod}}}

\State{$\mat{Q}^{(n)}{\leftarrow}\lambda \mat{W}^{(n)}+\sum_{p=1}^N\text{MatVec} (\mat{A}^{(n)}{,}\mat{A}^{(p)}{,}\mat{\Gamma}^{(n,p)}{,}\mat{W}^{(p)})$}
\EndFor

\State{$\alpha \leftarrow{\sum_{n=1}^N\langle\mat{R}^{(n)},\mat{Z}^{(n)}\rangle}/ {\sum_{n=1}^N\langle\mat{W}^{(n)}, \mat{Q}^{(n)}\rangle }$}

\For{\texttt{$n \in \inti{1}{N} $}}

\State{$\mat{V}^{(n)} \leftarrow \mat{V}^{(n)} + \alpha\mat{W}^{(n)}$}

\State{$\mat{R}^{(n)} \leftarrow \mat{R}^{(n)} - \alpha\mat{Q}^{(n)}$
}

\State{$\mat{Z}^{(n)}\leftarrow \mat{R}^{(n)}\mat{P}_{\text{inv}}^{(n)}$} 

\EndFor

\State{$\beta \leftarrow {\sum_{n=1}^N \langle \mat{R}^{(n)},{\mat{Z}^{(n)}}} \rangle/{\sum_{n=1}^N \langle \mat{W}^{(n)}, \mat{Q}^{(n)}\rangle} $ 
}

\For{\texttt{$n \in \inti{1}{N} $}}
\State{$\mat{W}^{(n)} \leftarrow \mat{Z}^{(n)}+ \beta \mat{W}^{(n)}$
}

\EndFor

\EndWhile

\State{\Return updates $\{ \mat{V}^{(1)}, \ldots , \mat{V}^{(N)}\}$ to factor matrices }

\end{algorithmic}

\end{algorithm}

\subsection{Preconditioning for Conjugate Gradient}

Preconditioning is often used to reduce the number of iterations in conjugate gradient.
For CP decomposition, the structure of the Gauss-Newton approximate Hessian $\mat{H}=\mat{J}^T\mat{J}$ admits a natural block-diagonal Kronecker product preconditioner~\cite{phan2013low}.
Each of the $N$ diagonal blocks $\mat{H}^{(n,n)}$ has a Kronecker product structure, $\mat{H}^{(n,n)}=\mat{\Gamma}^{(n,n)} \otimes \mat{I}$.
Consequently, its inverse is 
\[{\mat{H}^{(n,n)}}^{-1}={\mat{\Gamma}^{(n,n)}}^{-1} \otimes \mat{I},\]
which can be computed using $O(R^3)$ work per Gauss-Newton iteration and applied with $O(sR^2)$ cost per CG iteration.

We can also use the Cholesky factorization $\mat{\Gamma}^{(n,n)} = \mat{L} \mat{L}^T$,
\[\mat{H}^{(n,n)} = \mat{\Gamma}^{(n,n)} \otimes \mat{I} = (\mat{L} \mat{L}^T) \otimes \mat{I} = (\mat{L} \otimes \mat{I})(\mat{L}^T \otimes \mat{I}), 
\]
in which case application of ${\mat{H}^{(n,n)}}^{-1}$ can be applied in a stable way via triangular solve.
However, we found that performing triangular solves via ScaLAPACK~\cite{Dongarra:1997:SUG:265932} is a bottleneck for parallel execution (see the weak scaling results in Section~\ref{sec:exp}), as backward and forward substitution have polynomial depth.
Consequently, we compute the inverse of $\mat{\Gamma}^{(n,n)}$ and use tensor contractions to apply it in our parallel implementation.

\subsection{Complexity comparison between ALS and GN}

We present the cost of our Gauss-Newton implementation in Table~\ref{table:compare}. The right hand side of the Gauss-Newton iteration is the gradient of the residual function, which can be calculated using dimension trees similar to the ALS algorithm, thus requiring $O(s^NR)$ amount of work, the same as an ALS sweep.
With the use of implicit CG to solve the linear least squares problems in Gauss-Newton, the cost can be dominated by CG iterations, each of which requires $O(N^2sR^2)$ work. In exact arithmetic, CG should converge in at most $NsR$ iterations.

We compare this iterative approach to the best known methods for direct inversion of the approximate Hessian for CP decomposition with Gauss-Newton.
These approaches exploit the block structure of the approximate Hessian matrix, achieving a cost of $O(N^3R^6)$~\cite{phan2013low}, which may be improved to $O(NR^6)$ at the sacrifice of some numerical stability~\cite{tichavsky2013further}.
These methods accelerate inversion when $R$ is small.

However, CP decomposition may be accelerated by an initial Tucker factorization to decrease any dimension greater than $R$ down to $R$.
Tucker preserves exact CP rank and is easier to compute than CP (HoSVD is exact provided existence of an exact Tucker decomposition and is near-optimal for approximation).
When $R\geq s$, it is less clear whether the iterative or direct method is preferred.
One overhead of the direct approach is a memory footprint overhead of $O(NR^4)$.

\begin{table}[t]
\caption{Cost comparison between dimension tree ALS and Gauss-Newton methods.
Depth is quantified with $\tilde{O}$ to omit logarithmic depth factors associated with summations.}
\label{table:compare}
\centering
\begin{adjustbox}{width=0.7\textwidth}
\begin{tabular}{|c|c|c|}
  \hline
 \textbf{method} & \textbf{work} & \textbf{depth} \\
 ALS dimension tree~\cite{kaya2016parallel} & $O(s^NR + NR^3)$ & $\tilde{O}(N+R)$ \\
 GN with Cholesky~\cite{paatero1997weighted} & $O(s^NR + N^3s^3R^3)$ & $\tilde{O}(NsR)$ \\
 GN with fast inverse~\cite{phan2013low} & $O(s^NR + N^3R^6)$ & $\tilde{O}(NR^2)$ \\
 GN with faster inverse~\cite{tichavsky2013further} & $O(s^NR + NR^6)$ & $\tilde{O}(R^2)$ \\
 Implicit GN CG step~\cite{sorber2013optimization} & $O(N^2sR^2)$ & $\tilde{O}(1)$ \\
 GN step with $I$ CG iter  & $O(s^NR + IN^2sR^2)$ & $\tilde{O}(N+R+I)$ \\
 GN step with exact CG  & $O(s^NR + N^3s^2R^3)$ & $\tilde{O}(NsR)$ \\
  \hline
\end{tabular}
\end{adjustbox}
\end{table}

We quantify the work and depth (number of operation along critical path, lower bound on parallel cost)
of ALS and alternative methods for Gauss-Newton in Table~\ref{table:compare}.
The depth analysis for Gauss-Newton with CG assumes use of preconditioning with explicit inverse computation.
To quantify the depth of direct linear system solves (necessary in ALS and direct Gauss-Newton), we assume standard approaches (e.g., Gaussian elimination), which have a depth equal to matrix dimension, as opposed to polylogarithmic-depth matrix inversion methods~\cite{csanky1975fast}.
The communication costs associated with ALS and Gauss-Newton methods can be reduced to known analysis
for MTTKRP~\cite{ballard2018communication}, matrix multiplication~\cite{solomonik2011communication}, and Cholesky factorization~\cite{ballard2010communication}.
This analysis of cost and depth suggests that Gauss-Newton with implicit CG achieves the best cost and parallelism among Gauss-Newton variants when $s=O(R)$.
However, ALS generally offers more parallelism than Gauss-Newton with implicit CG, when the number of CG iterations is sufficiently large so as to dominate cost.


\section{Implementation}
\label{sec:implement}

We implement both dimension tree based ALS algorithm and Gauss-Newton algorithm in Python\footnote{Our implementations are publicly available at \url{https://github.com/cyclops-community/tensor_decomposition}.}. We leverage a backend wrapper for both NumPy and the Python version of Cyclops Tensor Framework~\cite{solomonik2014massively}, so that our code can be tested and efficiently executed both sequentially and with distributed-memory parallelism for tensor operations.
In addition, we write both the ALS and Gauss-Newton optimization algorithms in an optimizer class, and each ALS sweep / Gauss-Newton iteration is encapsulated as a step member function in the optimizer class. This framework can be easily extended to included other optimization algorithms for tensor decompositions.
Cyclops provides a high-level abstraction for distributed-memory tensors, including arbitrary tensor contractions and matrix factorizations such as Cholesky and SVD via ScaLAPACK~\cite{Dongarra:1997:SUG:265932}.
The ALS implementation is based on previous work~\cite{ma2018accelerating} and uses dimension trees to minimize cost.

Our tensor contraction formulation of the Gauss-Newton method makes it easy to implement with NumPy and Cyclops.
Both libraries provide an \lstinline[language=Python]{einsum} routine for tensor contractions specified in Einstein summation notation.
Using this routine, the Gauss-Newton method can be specified succinctly as in the following code snippet, where lists of tensors are used to store the factor matrices $\mat{A}^{(n)}$, components of the input and output matrices (set of vectors) $\mat{W}^{(p)}$ and $\mat{U}^{(n)}$, and matrices $\mat{\Gamma}^{(n,p)}$.
\begin{lstlisting}[language=Python, caption={Implicit Matrix-Vector Product in GN Method}, captionpos=b, label={lst:gn}]
u = []
for n in range(N):
  u.append(zeros((s,R)))
  for p in range(N):
    if n == p:
      U[n] += einsum("rz,kz->kr",Gamma[n,p],W[p])
    else:
      U[n] += einsum("kz,lr,rz,lz->kr", \
                     A[n],A[p],Gamma[n,p],W[p])
\end{lstlisting}

Our current implementation does not parallelize over the $N^2$ matrix vector products. However, for the case of equidimensional tensors, we can cast the list of factor matrices as a tensor and cast the above contractions into two tensor contractions to achieve parallelization over the $N^2$ contractions. In the following code snippet we have the batched contraction where the input and output list of matrices are tensors $\tsr{V}$ and $\tsr{R}$ respectively, $\tsr{D}$ is the list of $\mat{\Gamma}^{(n,n)}$, $\tsr{G}$ is a fourth order tensor where the entries along the first two modes, $n^{th}$ and $p^{th}$ mode is $\mat{\Gamma}^{(n,p)}$ when $n \neq p$ and a matrix of zeros of size $R \times R$ along the diagonal and $\tsr{A}$ is the list of factor matrices:

\pagebreak

\begin{lstlisting}[language=Python, caption={Implicit Matrix-Vector Product with batched tensor contractions}, captionpos=b, label={lst:gn2}]
R = einsum("niz,nzr->nir",V,D)
R += einsum("niz,pjr,npzr,pjz->nir",A,A,G,V)

\end{lstlisting}

 Note that in our implementation, extra work is spent due to the zeros on the diagonal of $\tsr{G}$ tensor and computing the contractions with diagonal terms of the Hessian sequentially. The performance can be further improved by concurrent contraction of these terms.

\section{Numerical Experiments}
\label{sec:exp}

We performed numerical experiments to compare the performance of dimension tree based ALS algorithm and Gauss-Newton algorithm on both synthetic and application tensors. 
Our experiments consider
four types of tensors:

\noindent \textbf{Tensors made by random matrices}. We create these tensors based on known uniformly distributed randomly-generated factor matrices $\mat{A}^{(n)}\in (a,b)^{s\times R}$, 
   \(
    \tsr{X} = [\![ \mat{A}^{(1)}, \ldots , \mat{A}^{(N)} ]\!]. 
    \)

\noindent \textbf{Tensors made by Gaussian matrices}. We create tensors based on known Standard Gaussian distributed randomly-generated factor matrices $\mat{A}^{(n)}\in \mathcal{N}(0,1)^{s\times R}$, 
   \(
    \tsr{X} = [\![ \mat{A}^{(1)}, \ldots , \mat{A}^{(N)} ]\!]. 
    \)
    
\noindent  \textbf{Quantum chemistry tensors}.
    We also consider the density fitting tensor (Cholesky factor of the two-electron integral tensor) arising in quantum chemistry.
    Its CP decomposition yields the tensor hypercontraction format of the two-electron integral tensor, which enables reduced computational complexity for a number of post-Hartree-Fock methods~\cite{hohenstein:221101}.
    Acceleration of CP decomposition for this quantity has previously been a subject of study in quantum chemistry~\cite{hummel2017low}.
    We leverage the PySCF library~\cite{sun2018pyscf} to generate the three dimensional compressed density fitting tensor, representing the compressed restricted Hartree-Fock wave function of a water molecule chain systems with a
STO-3G
basis set. We vary the number of molecules in the system from 3 to 40, comparing the efficacy of ALS and Gauss-Newton method under different settings.

\noindent \textbf{Matrix multiplication tensor}.
    A hard case for CP decomposition is the matrix multiplication tensor, defined as an order three unfolding (combining pairs of consecutive modes) of
    \[t_{ijklmn} = \delta_{lm}\delta_{ik}\delta_{nj}.\]
    This tensor simulates multiplication of matrices $\mat{A}$ and $\mat{B}$ via
    \[c_{ij}=\sum_{klmn} t_{ijklmn}a_{kl}b_{mn} = \sum_{l} a_{il}b_{lj}.\]
    Its exact CP decompositions give different bilinear algorithms for matrix multiplication, including classical matrix multiplication with rank $s^{3/2}$ and Strassen's algorithm~\cite{Strassen_1969} with rank $s^{\log_4(7)}$.
    Determining the minimal CP rank for multiplication of $n$-by-$n$ matrices with $n\geq 3$ (so $s\geq 9$) is an open problem~\cite{Pan:1984:MMF:2212} that is of interest in theory and practice.

To maintain consistency throughout the experiments, we run CG until a relative tolerance of $10^{-3}$. We use the metrics \textit{relative residual} and \textit{fitness} to evaluate the convergence. Let $\Tilde{\tsr{X}}$ denote the tensor reconstructed by the factor matrices, the relative residual and fitness are 
\begin{align*}
    r = \frac{\|\tsr{X} - \Tilde{\tsr{X}}\|_F}{\|\tsr{X}\|_F}, \quad
    f = 1 - \frac{\|\tsr{X} - \Tilde{\tsr{X}}\|_F}{\|\tsr{X}\|_F}.
\end{align*}
We collect our experimental results with NumPy backend on the Blue Waters Supercomputer, and our results with Cyclops backend on the Stampede2 supercomputer of Texas
Advanced Computing Center located at the University of Texas at Austin using XSEDE~\cite{6866038}. 

On Blue Waters, we use one processor of the XE6 dual-socket nodes for each sequential experiment with NumPy backend.  
On Stampede2, we leverage the Knigt’s Landing (KNL) nodes exclusively, each of which consists of 68 cores,
96 GB of DDR RAM, and 16 GB of MCDRAM. These nodes are connected via a 100
Gb/sec fat-tree Omni-Path interconnect. We use Intel compilers and the MKL library
for BLAS and batched BLAS routines within Cyclops.
We use 64 processes per node on Stampede2 for all experiments. 

We study the effectiveness of ALS and Gauss-Newton on CP decomposition based on the following metrics:

\noindent \textbf{Convergence likelihood.}
    We compare the likelihood of the CP decomposition to recover the original low rank structure of the input tensor with both algorithms. 

\noindent \textbf{Convergence behavior.}
    We compare the convergence progress w.r.t.\ execution time of ALS and Gauss-Newton for all the tensors listed above. Experiments are performed with NumPy backend for small and medium-sized tensors, while the Cyclops backend is used for large tensors.

\noindent \textbf{Parallel Performance.} We perform a parallel scaling analysis to compare the simulation time for one ALS sweep of the dimension tree based ALS algorithm and the conjugate gradient iteration of the Gauss-Newton algorithm.

\begin{figure*}[t]
\centering
\begin{subfigure}[Tensor made by random factor matrices by elements in $(0,1)$] {\label{fig:proba}\includegraphics[width=0.45\textwidth, keepaspectratio]{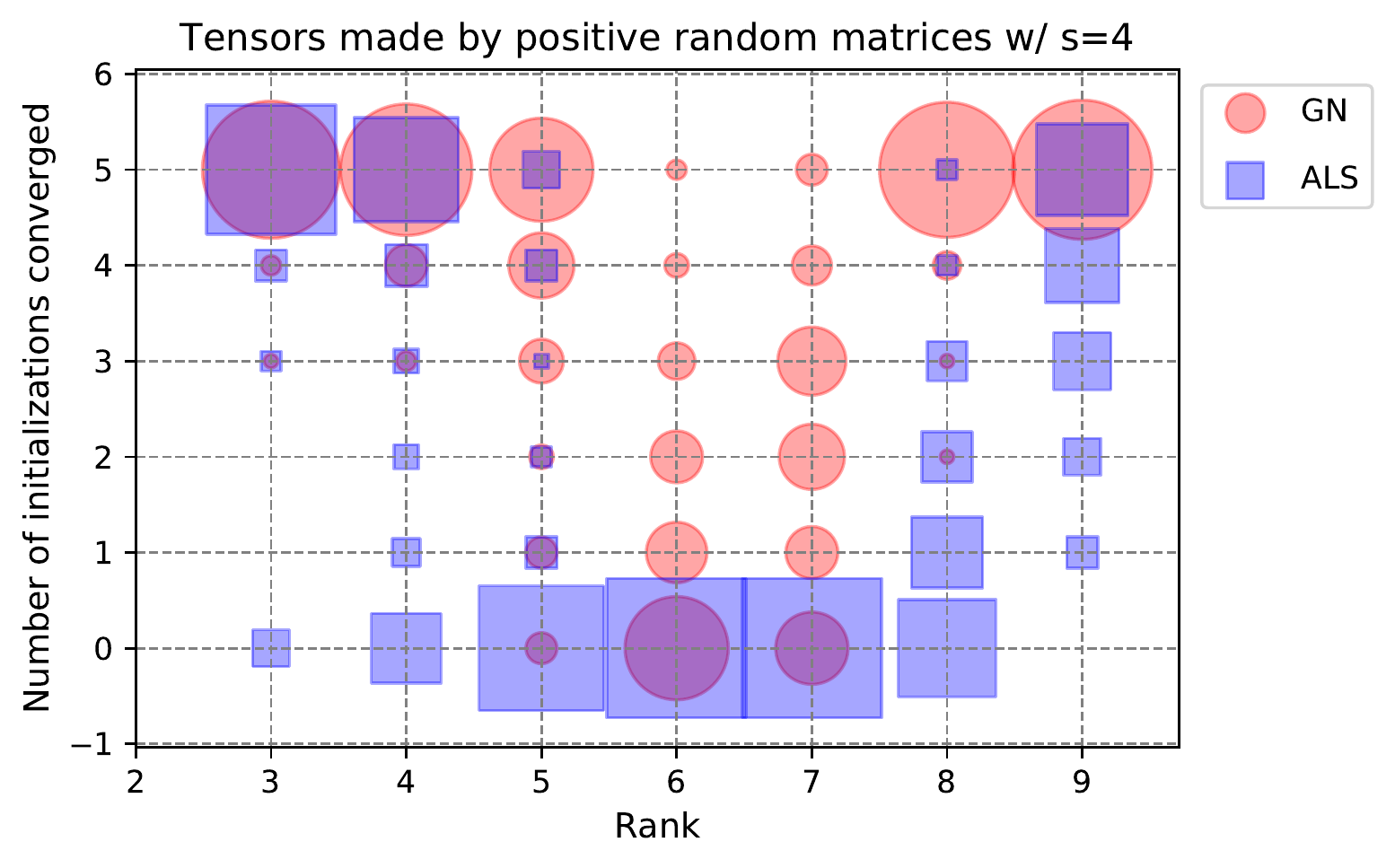}} \end{subfigure}
\begin{subfigure}[Tensor made by random factor matrices by elements in $(0,1)$] {\label{fig:resplot}\includegraphics[width=0.45\textwidth, keepaspectratio]{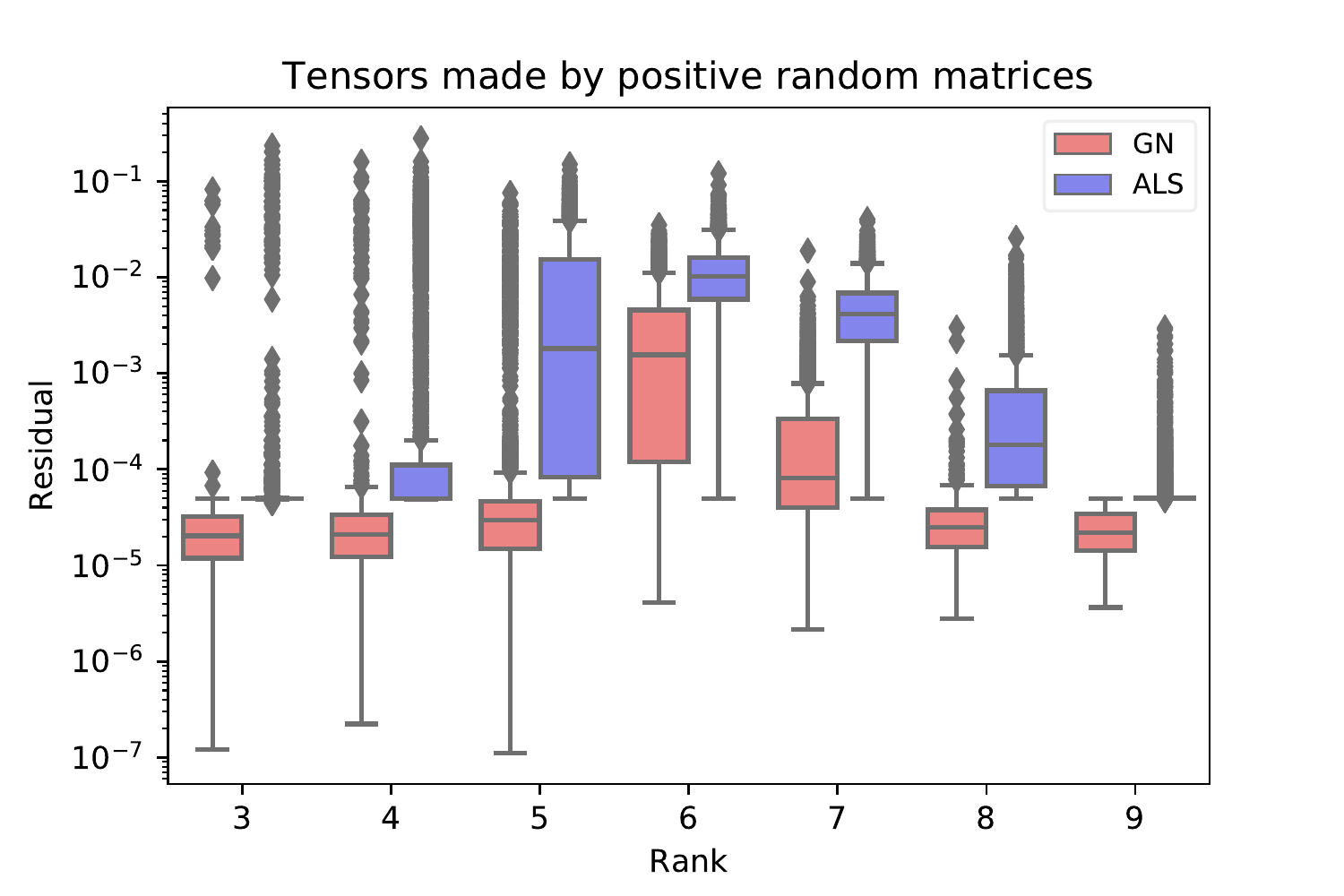}} \end{subfigure}
\caption{Convergence of Gauss-Newton with varying identity regularization and ALS algorithms for recovery of exact CP decomposition with random positive factor matrices.}
\label{fig:posi_area}
\end{figure*}

\begin{figure*}[t]
\centering
\begin{subfigure}[Convergence results for random tensor with 5 and 15 initializations] {\label{fig:negprob}\includegraphics[width=0.45\textwidth, keepaspectratio]{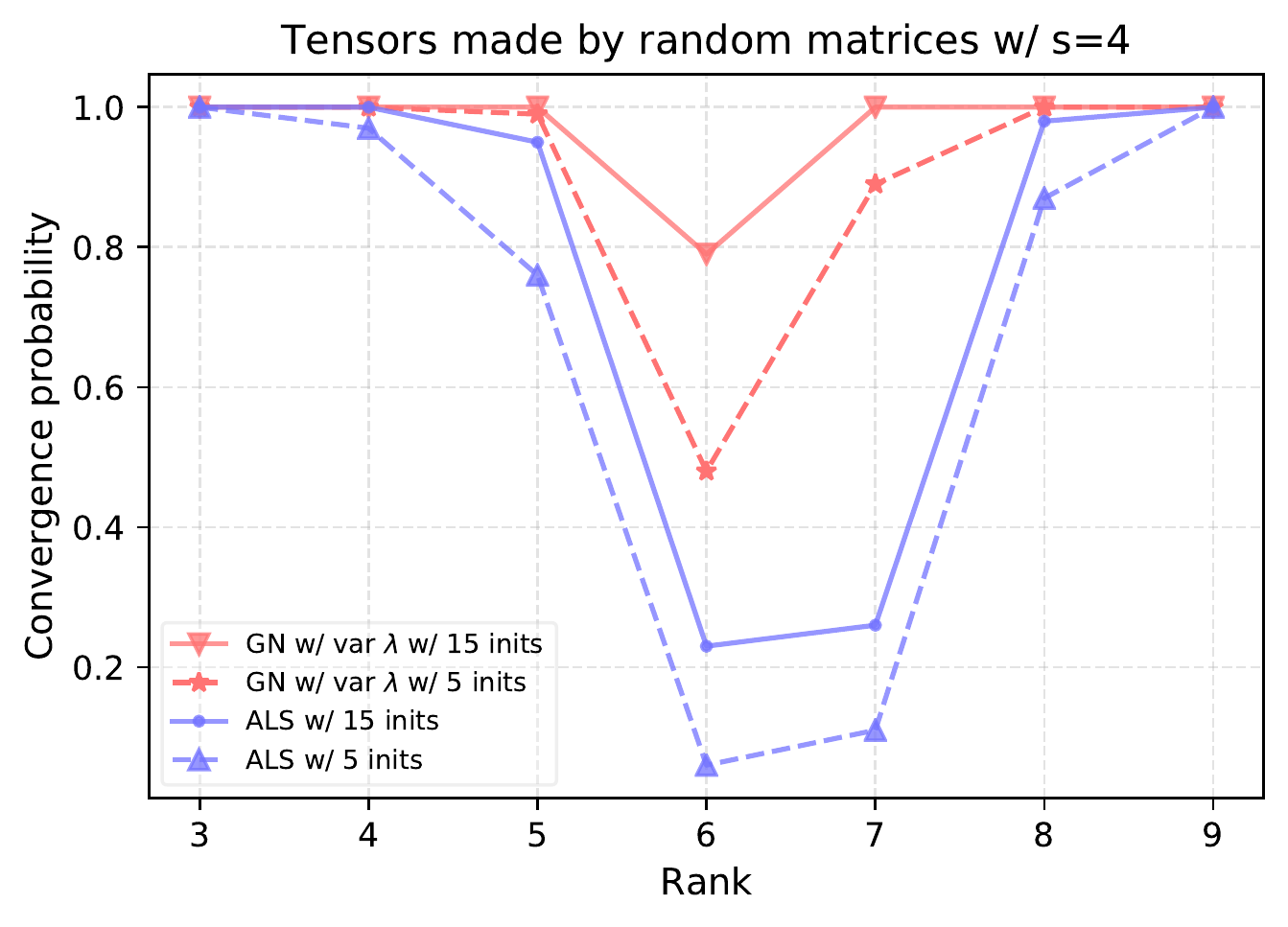}} \end{subfigure}
\begin{subfigure}[Convergence results of all the variants for Gaussian random tensor] {\label{fig:allgauss}\includegraphics[width=0.45\textwidth, keepaspectratio]{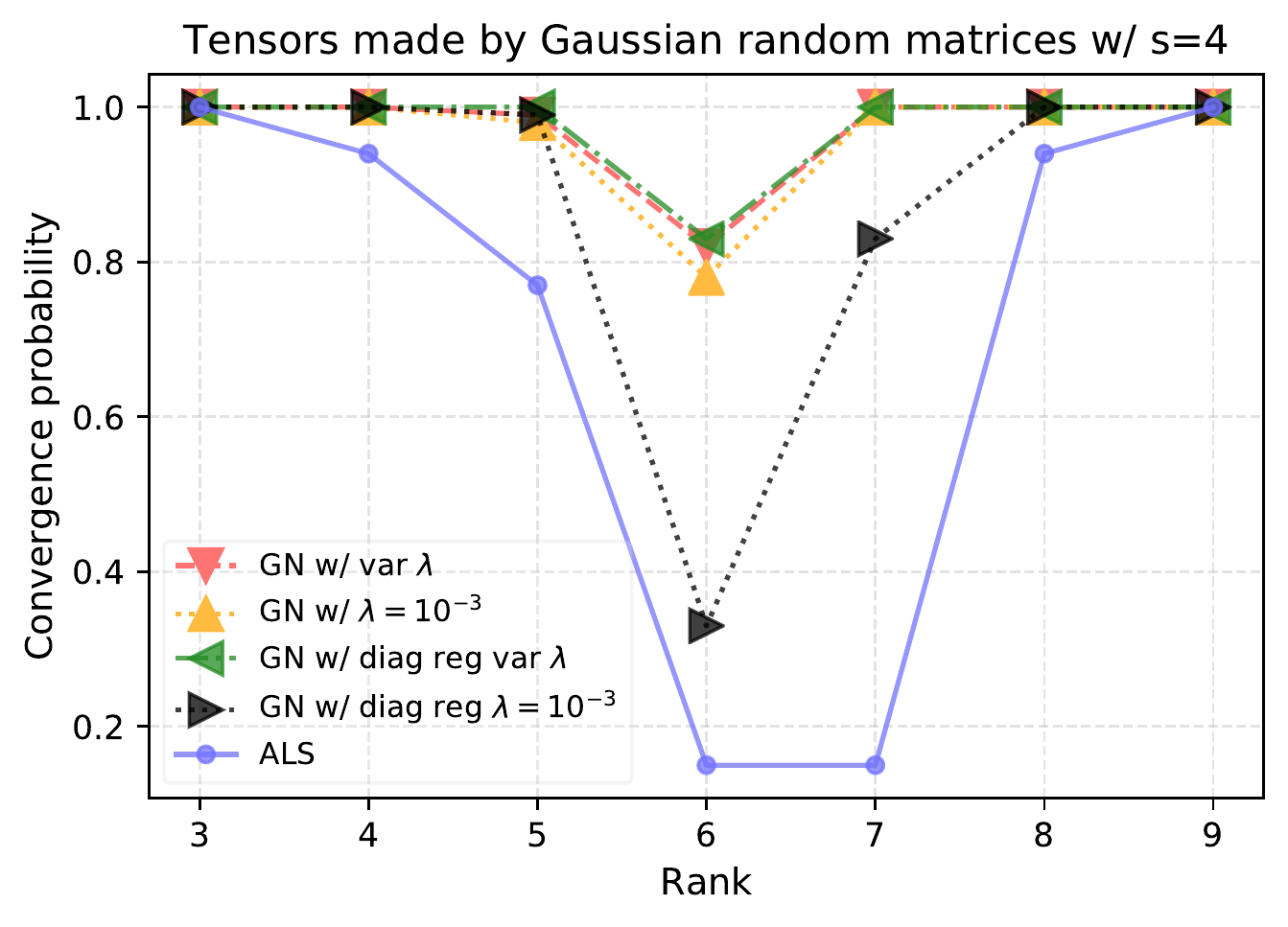}} \end{subfigure}
\begin{subfigure}[Convergence results of all the variants for positive random tensor] {\label{fig:allrand}\includegraphics[width=0.45\textwidth, keepaspectratio]{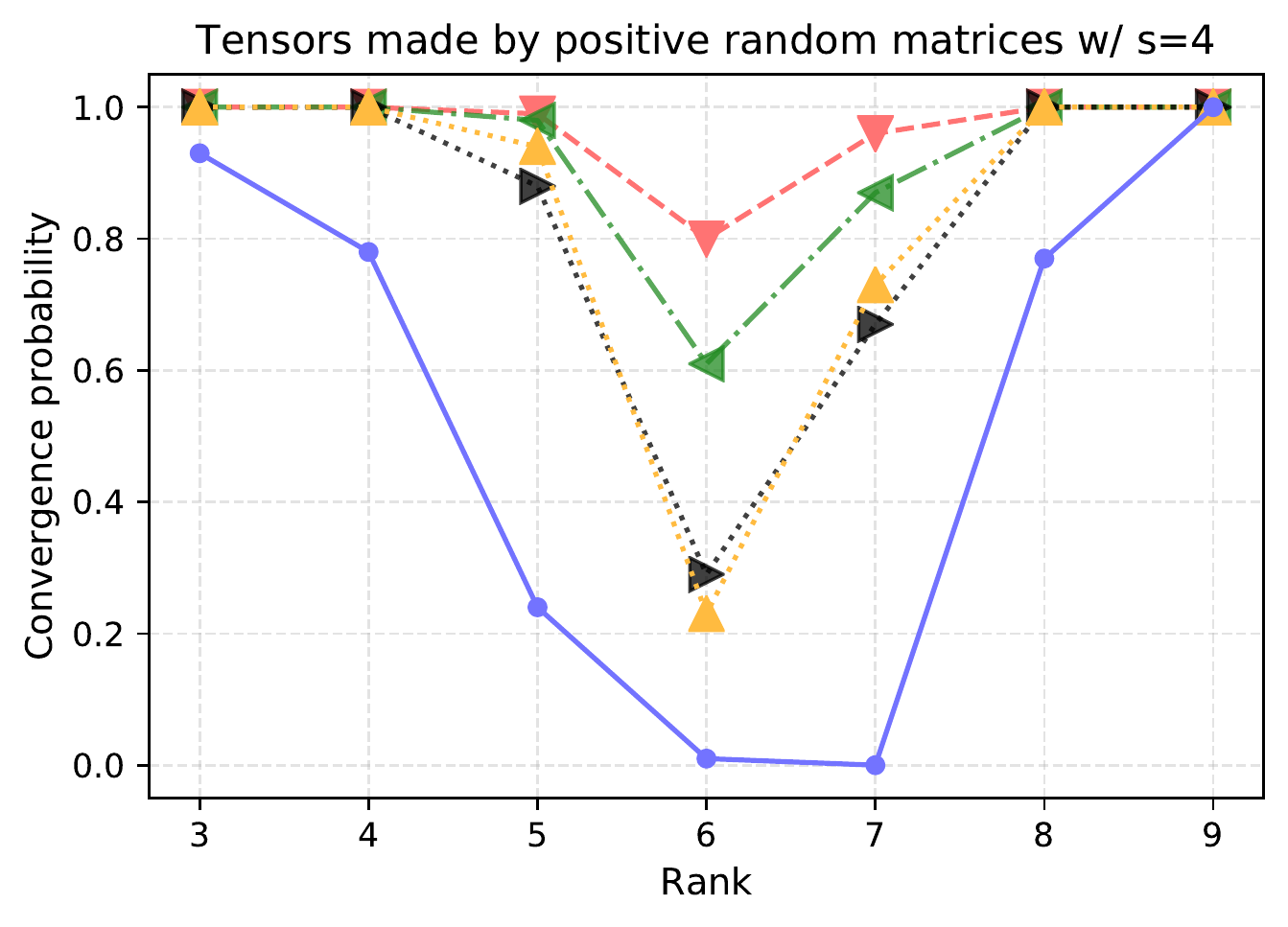}} \end{subfigure}
\begin{subfigure}[Convergence results of all the variants for random tensor] {\label{fig:allnegrand}\includegraphics[width=0.45\textwidth, keepaspectratio]{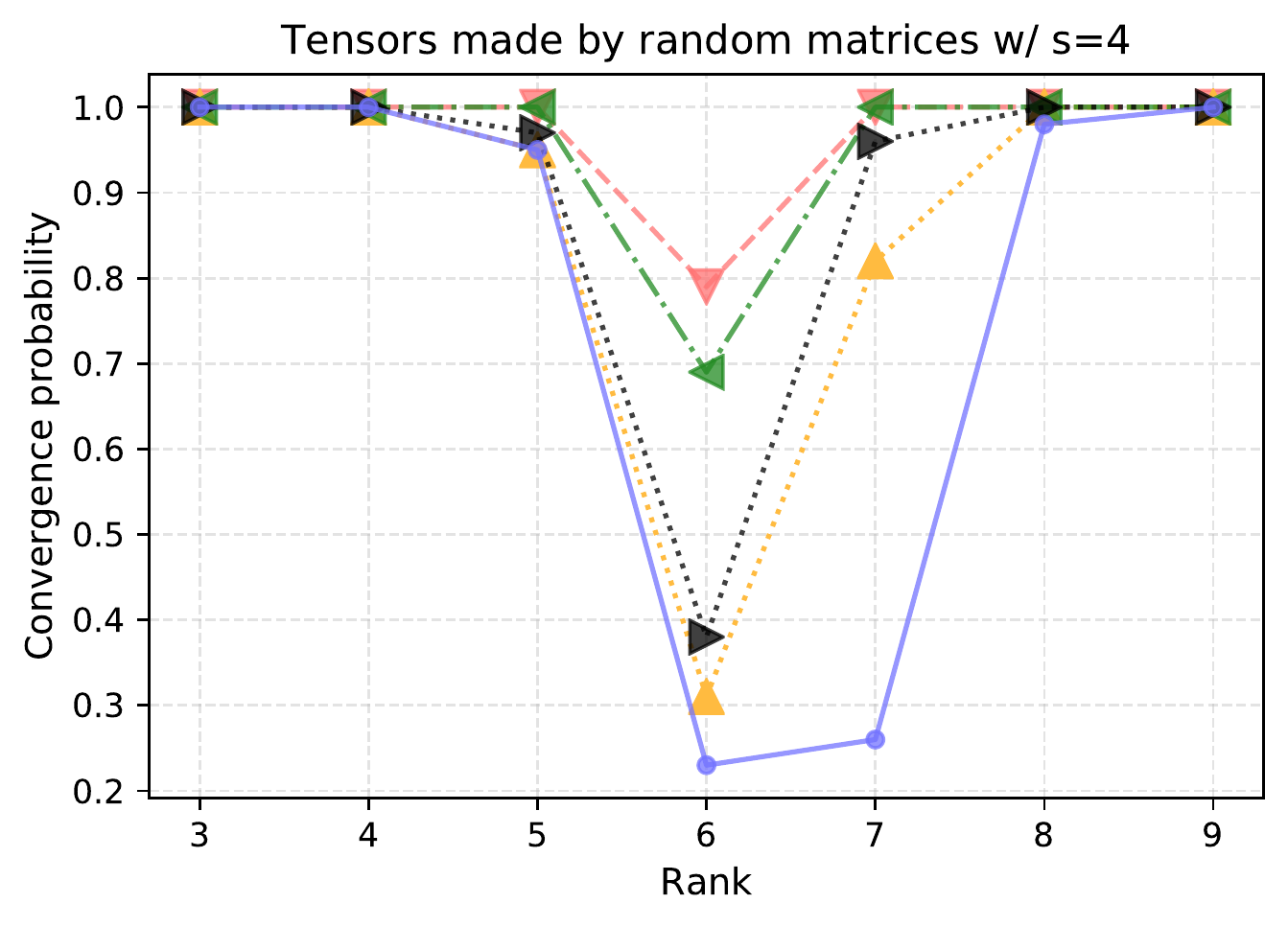}} \end{subfigure}
\caption{Convergence results of various versions of regularization of Gauss-Newton and ALS for recovery of exact CP decomposition with
factor matrices with entries selected using uniform random, positive uniform random, and Gaussian distributions.
}
\label{fig:prob_plots}
\end{figure*}

\begin{figure*}[t]
\centering
\begin{subfigure}[Tensors made by factor matrices with standard Gaussian distribution] {\label{fig:Regplot}\includegraphics[width=0.55\textwidth, keepaspectratio]{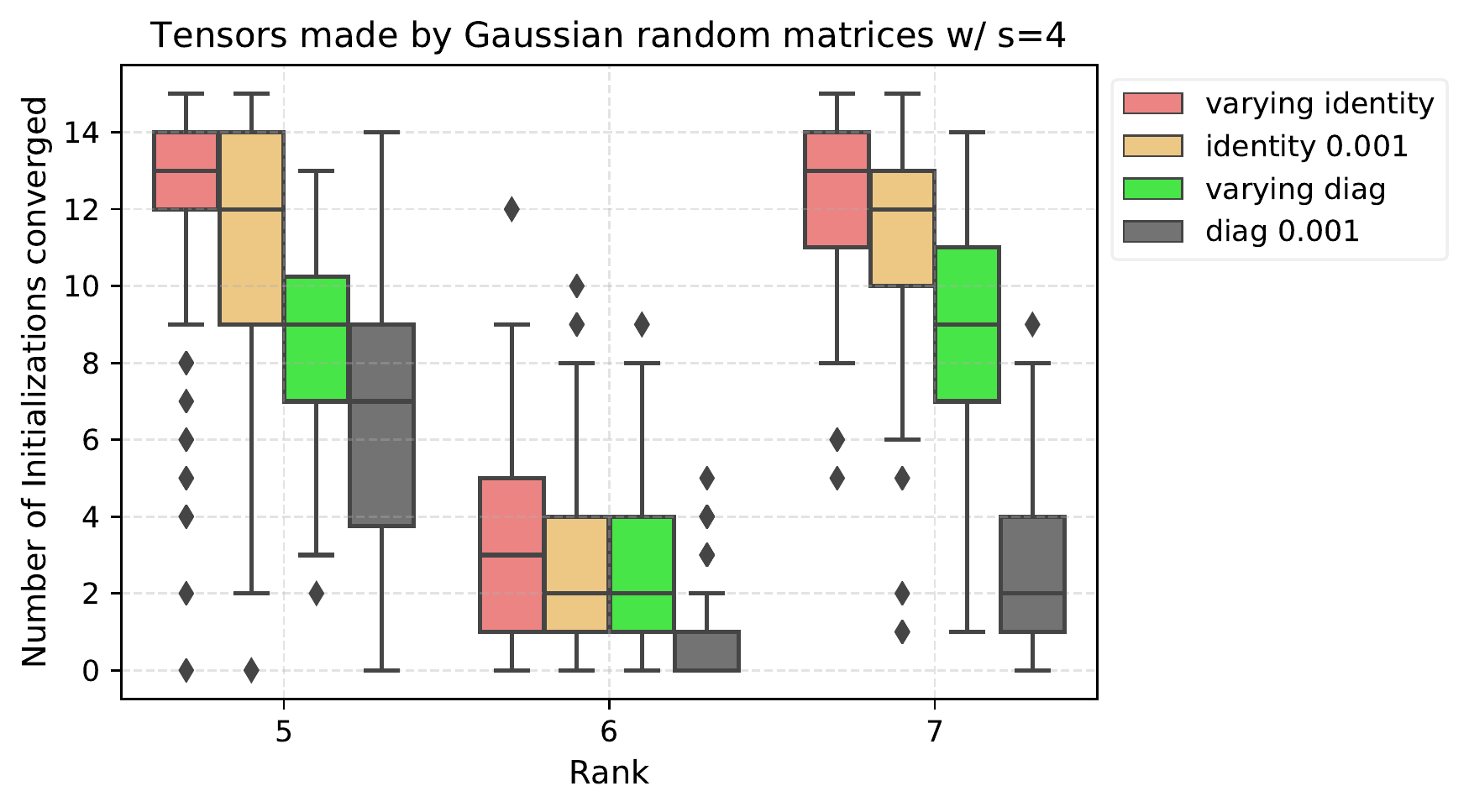}} \end{subfigure}
\begin{subfigure}[Matrix multiplication tensors] {\label{fig:matmul}\includegraphics[width=0.43\textwidth, keepaspectratio]{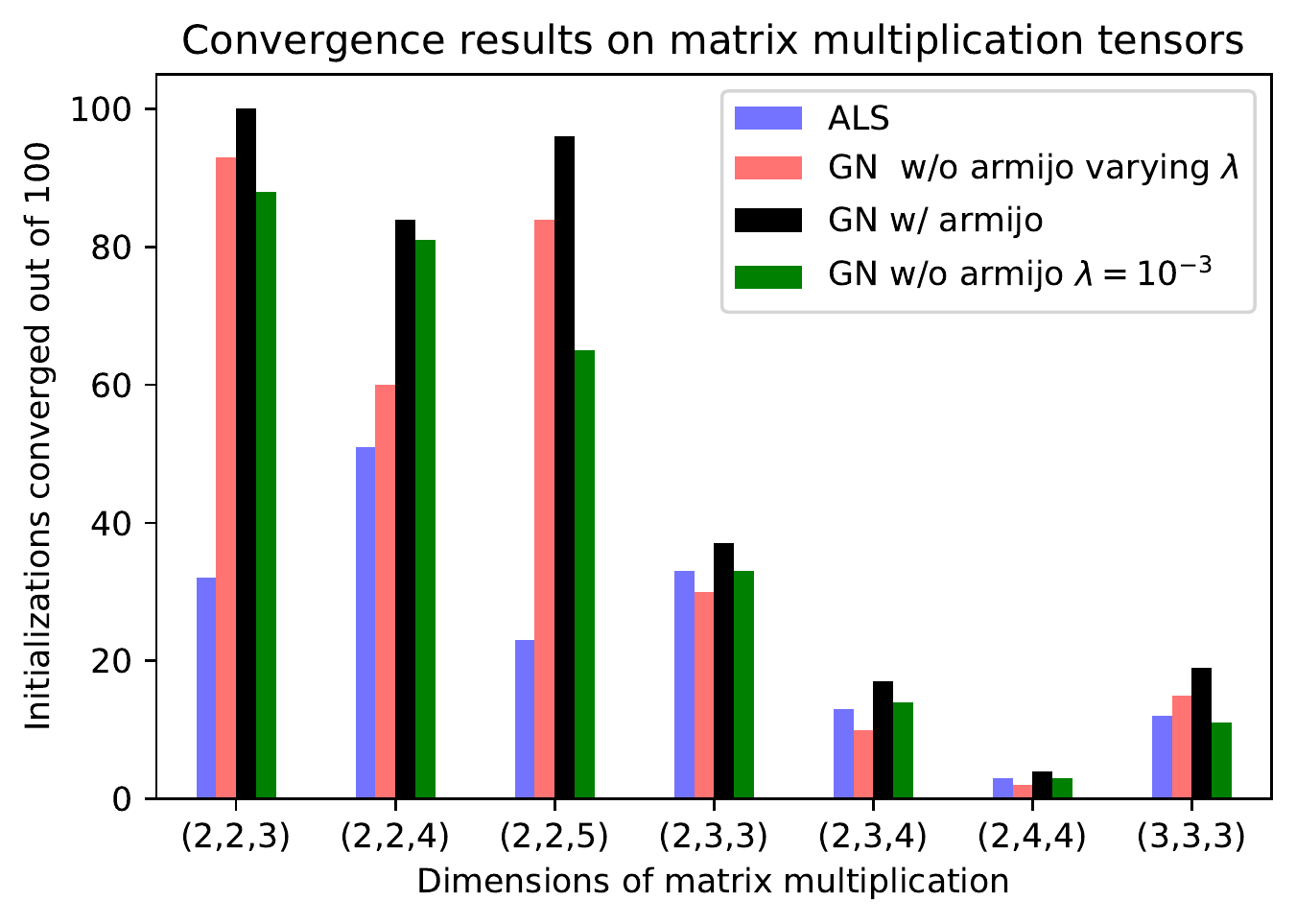}} \end{subfigure}
\caption{Convergence results of various versions of regularization of Gauss-Newton and ALS for recovery of exact CP decomposition for Gaussian random tensors and matrix mutliplication tensors.}
\label{fig:matmul}
\end{figure*}

\subsection
{Convergence likelihood}
\label{subsec:convlikelihood}

We compare the convergence likelihood of CP decomposition for random low-rank tensors, optimized with ALS algorithm and Gauss-Newton algorithm with constant and varying regularization. We run the algorithms until the residual norm is less than $5 \times 10^{-5}$, or the norm of the update is less than $10^{-7}$,  or a maximum of 500 and 10,000 iterations for Gauss-Newton and ALS, respectively for 100 problems.
The results are presented in Figure~\ref{fig:posi_area},~\ref{fig:prob_plots} and~\ref{fig:matmul}. 
We set the tensor order $N=3$, size in each dimension $s=4$, and compare the convergence likelihood under different CP ranks.
These results are representative of behavior observed across a variety of choices of $s$ and $R$.

In Figure~\ref{fig:proba} and ~\ref{fig:resplot}, we run Gauss-Newton and ALS with factor matrices sampled from $(0,1)$ uniformly at random with 5 initializations each for CP rank ranging from 3 to 9. 
The diameter of the circle and the
side length of the square are proportional to the number of problems converged for the corresponding number of initializations in ~\ref{fig:proba}. It is evident that Gauss-Newton exhibits a higher probability of convergence than ALS as the circles are always bigger than the squares for higher number of initializations converged. We can observe in Figure~\ref{fig:resplot} that Gauss-Newton with varying regularization is more likely to reach a lower residual when compared with ALS (giving both ample number of iterations).

In Figure~\ref{fig:negprob}, we run both algorithms with factor matrices sampled from $(-1,1)$ uniformly at random with 5 and 15 initializations. A point in the graph represents the probability of at least one initialization converging out of the total initializations. We observe similar behaviour over the various ranks, 6 being the most difficult to converge. However, with increasing the initializations we observe increase in the convergence probability for both the algorithms but Gauss-Newton with identity varying regularization outperforms ALS.

In Figure \ref{fig:allgauss},~\ref{fig:allrand} and~\ref{fig:allnegrand}, we compare all the algorithms with different types of tensors with 15 initializations. These plots indicate that varying regularization improves convergence for both the variants of regularization in various types of tensors whereas ALS does not do well at convergence for these 2 types of tensors for the `harder' cases. Moreover, the probability for the varying identity regularization for remains constant various tensors suggesting that the convergence probability of the method is invariant to how the tensors were constructed.

In Figure~\ref{fig:Regplot}, we compare Gauss-Newton with different regularization techniques for tensors with factor matrices sampled from the standard Gaussian distribution for the `harder' cases (ranks 5 to 7) with $15$ initializations. Plotting the number of converged initializations per problem for these variants over the `harder' cases,
we observe that Gauss-Newton with varying identity regularization performs better than other variants of regularization. We also observe that varying the regularization parameter increases the number of converged problems which corroborates our claim that varying regularization improves the probability of convergence.

In Figure~\ref{fig:matmul}, we find the CP decomposition of matrix multiplication tensors with best known ranks~\cite{benson2015framework} with 100 initializations. For ALS algorithm, we start with a high regularization parameter, $\lambda = 0.01$ and decrease it gradually, by a factor of 2 after every 100 iterations, which is suggested in~\cite{smirnov2013bilinear} to increase the probability for finding the CP decomposition. We run ALS for 20000 iterations and the convergence criteria is set at $10^{-8}$. For Gauss-Newton method, we initialize it with 200 iterations of ALS with $\lambda = 0.01$ and then use Gauss-Newton with proposed regularization and with constant $\lambda = 10^{-3}$. We found that in this case using Armijo's condition~\cite{armijo1966}  for step-size control increases the probability of convergence for Gauss-Newton method which is more than both the constant and variable regularization strategy. We do not observe the same pattern where varying the regularization increases the convergence probability for Gauss-Newton in this case as the regularization is not on the $L^2$ norm of the factor matrices as opposed to ALS (which is shown to work better in this case). However, with Armijo's condition incorporated, convergence probability of Gauss-Newton is more than ALS with the mentioned regularization strategy.

\subsection{Parallel Performance}
\label{subsec:performance}

\begin{figure*}[t]
\centering
\begin{subfigure}[Weak scaling with fixed tensor size to number of processors ratio and tensor dimension to rank ratio] {\label{fig:bencha}\includegraphics[width=0.48\textwidth, keepaspectratio]{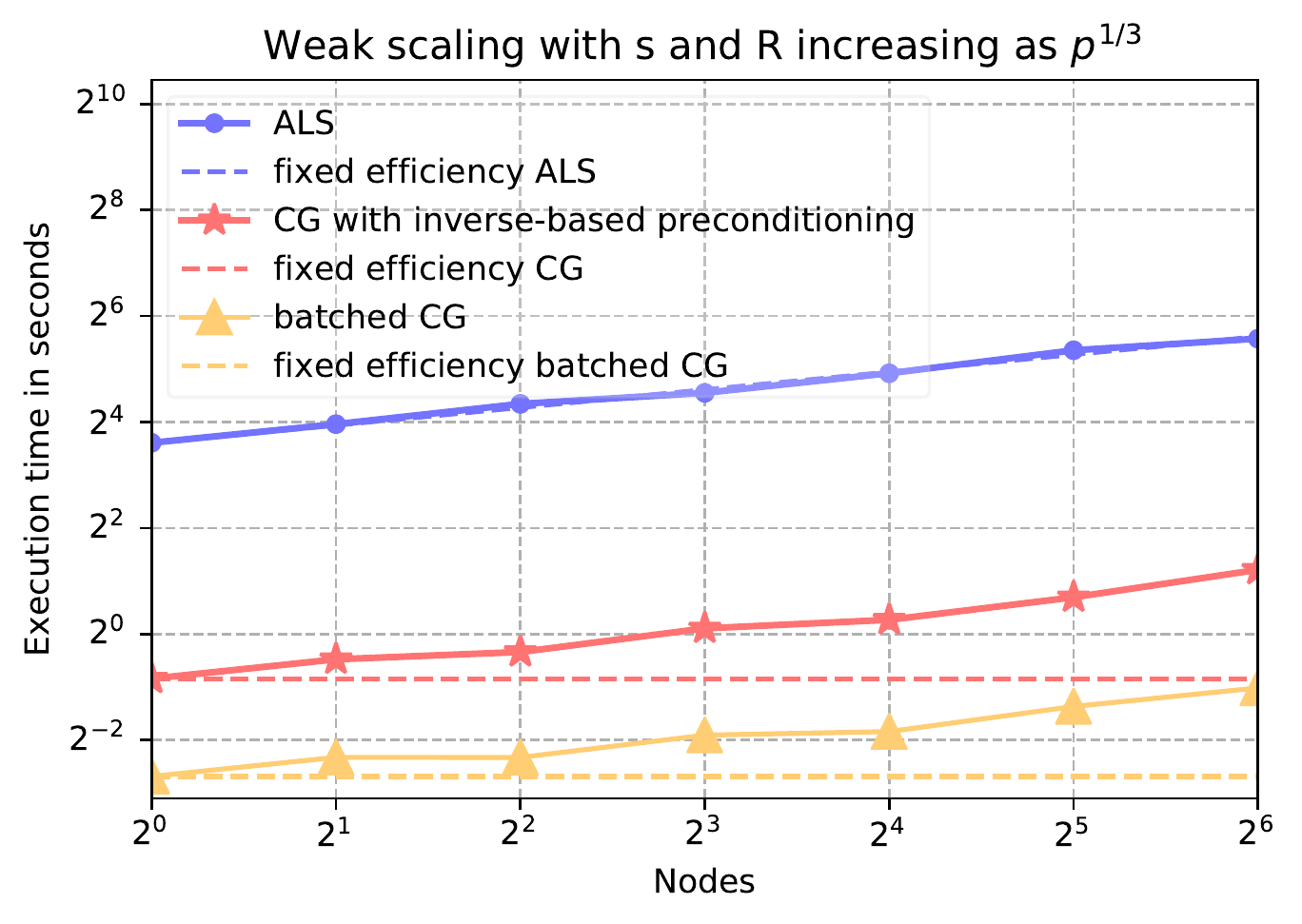}} \end{subfigure}
\begin{subfigure}[Weak scaling with fixed tensor size to processors ratio and compression ratio] {\label{fig:benchc}\includegraphics[width=0.48\textwidth, keepaspectratio]{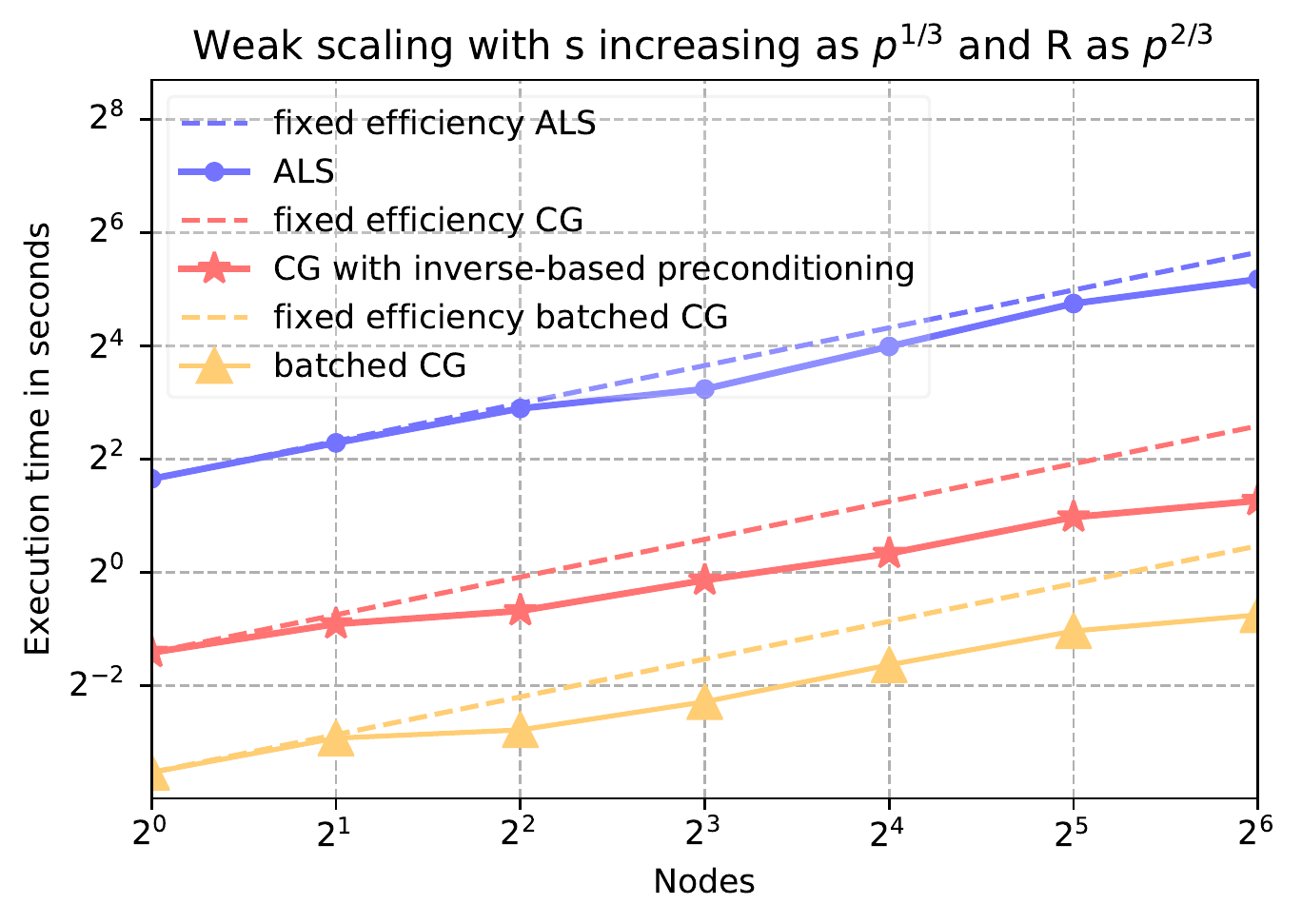}} \end{subfigure}

\begin{subfigure}[Strong scaling with fixed tensor size ] {\label{fig:benchb}\includegraphics[width=0.48\textwidth, keepaspectratio]{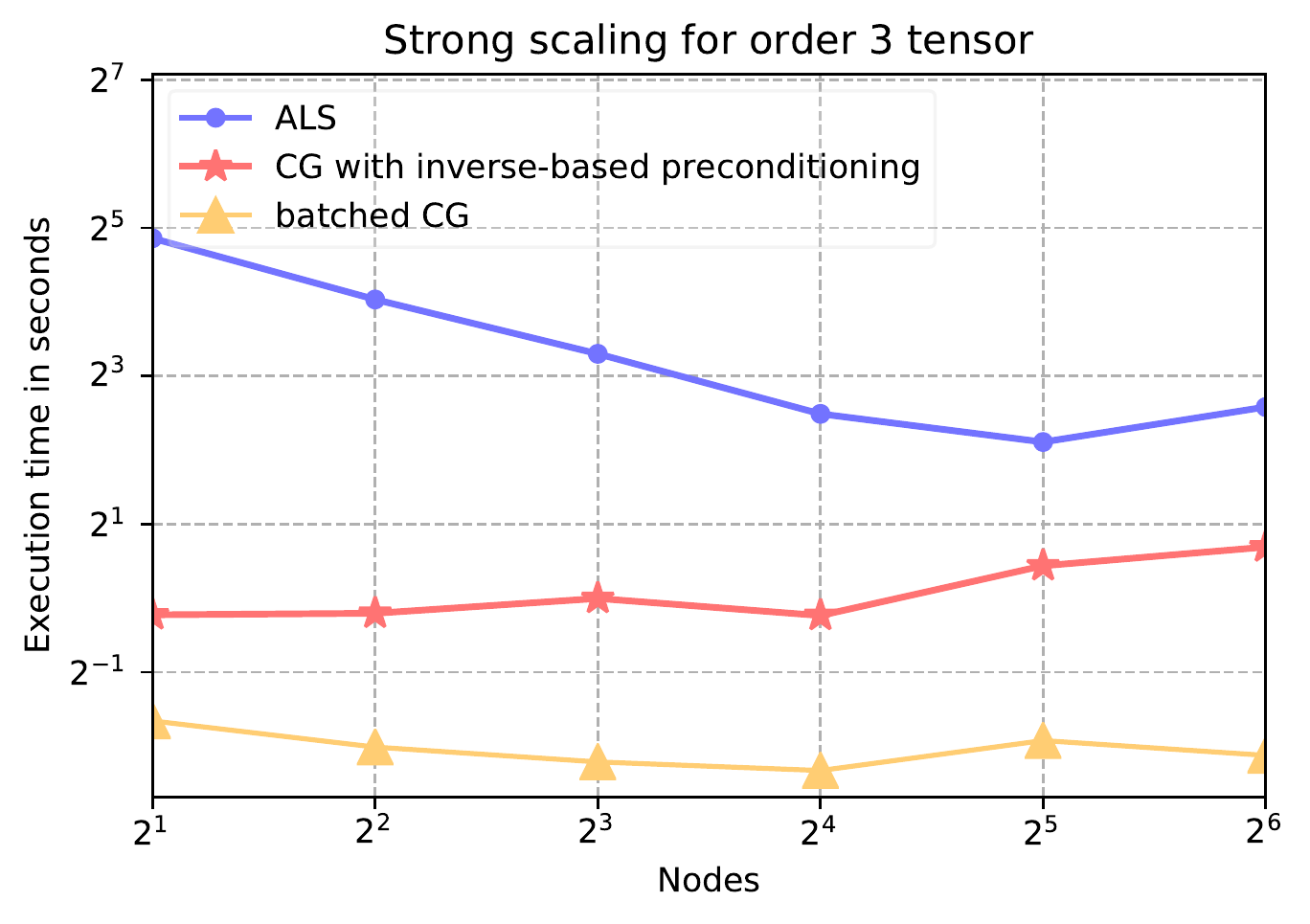}} \end{subfigure}

\caption{Benchmark results for one ALS sweep vs one CG iteration. Each data point is the mean time across 5 iterations.
}
\label{fig:bench}

\end{figure*}

We perform a parallel scaling analysis to compare the simulation time for one dimension tree based ALS sweep  and one conjugate gradient iteration of the Gauss-Newton algorithm. 
In Figure~\ref{fig:bencha}, we perform weak scaling with $p$ processors, considering order $N = 3$ tensors starting with dimension $s = 800$ and rank $R = 800$ and growing both as $p^{1/3}$ with increasing number of nodes $p$. This scaling maintains a fixed memory footprint of the tensor per processor, The work per processor for ALS being $O(s^3R/p)$, increases by a factor of $O(p^{1/3})$ with $p$ processors. For a CG iteration, the work per processor is $O(N^2sR^2/p)$, which remains constant per processor. Figure~\ref{fig:bencha}  shows that with the increase of number of nodes, the time for one ALS sweep scales perfectly while the efficiency for one conjugate gradient iteration drops to $24\%$ at 64 nodes due to the limited number of operations involved in the Hessian contraction. Also, note that the Hessian contraction takes up about half the time of a CG iteration as inner product, norm calculation take up a significant amount of time.
One conjugate gradient iteration is consistently around 20 times faster than one ALS sweep for different simulation sizes.
We observed that explicit calculation and use of the inverse eliminates a significant overhead compared to preconditioning using Cholesky and triangular solves.

We also implement 
batch CG which uses batched Hessian contractions as described in section~\ref{sec:implement}. Due to the fact that we extract more parallelism over the $N^2$ contractions by batching the contractions into a bigger tensor contraction, one CG iteration speeds up by a factor of 3.58 with 1 node and 4.68 with 64 nodes with this implementation.

In Figure~\ref{fig:benchc}, we perform weak scaling with $p$ processors, considering order $N = 3$ tensors starting with dimension $s = 600$ and rank $R = 300$, growing $s$ as $p^{1/3}$ and $R$ as $p^{2/3}$ with increasing number of nodes $p$. This scaling maintains a fixed memory footprint of the tensor and factor matrices per processor while the work per processor for ALS and CG  increases by a factor of $O(p^{2/3})$ per processor. Figure~\ref{fig:benchc}  shows that with the increase of number of nodes, the time for one ALS sweep and one conjugate gradient iteration increases and the efficiency improves with growing size. Although, Hessian contraction takes up only half of the time of a CG iteration, it takes $0.15$ seconds with 1 node and scales up with an efficiency of more than $200\%$, taking $1.12$ seconds with 64 nodes. The increase in efficiency is because of the increase in arithmetic intensity is increasing by $O(p^{2/3})$ per process, which leads to a speedup of greater than $O(p)$ for both the algorithms. These observations demonstrate a good weak scaling of CG iteration over increasing ranks.

For strong scaling, we consider order $N = 3$ tensors with dimension size $s = 1200$ and rank $R = 1200$.
 Figure~\ref{fig:benchb} shows that the conjugate gradient iteration time increases with the number of nodes, while the ALS sweep time decreases at first, and increases with more than 32 nodes due to  communication cost dominating afterwards. 
The CG iteration involves smaller matrix multiplications, and the contraction time does not scale with increasing node counts on account of the communication cost.
The Hessian contraction takes 0.45 seconds with 2 nodes and scales to 0.35 seconds with 16 nodes, while the operations including norm  calculation scales worsen as they are latency bounded and hence the time increases for the iteration. 
The batch CG contain bigger contractions, which results in improving the time and scaling as the batched contraction takes 0.22 seconds with 2 nodes and scales to 0.13 seconds with 16 nodes. The time taken here is also dominated by the norm and inner product calculations.
The ALS sweep is dominated by the MTTKRP calculations, which are more easily parallelizable and therefore make ALS achieve better scaling.
Overall, we observe that the Gauss-Newton CG iterations contain less parallelism than MTTKRP, but are weakly scalable when rank is increasing.

\subsection{Exact CP decomposition}
\label{subsec:accurate}

\begin{figure*}[t]
\centering
\begin{subfigure}[Tensor made by random factor matrices with elements distributed with standard normal distribution] {\label{fig:gauss80}\includegraphics[width=0.48\textwidth, keepaspectratio]{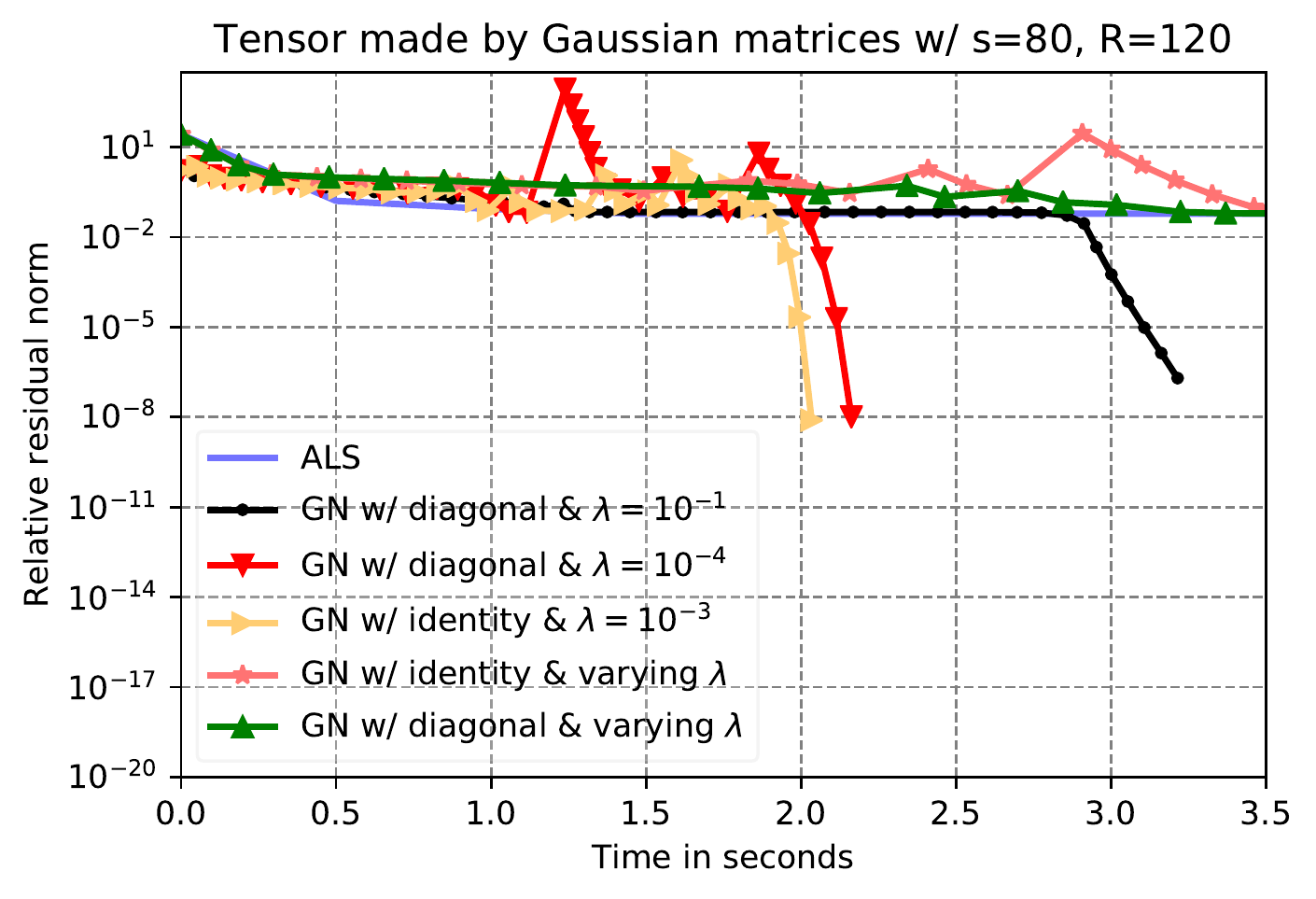}}\hspace{2mm} \end{subfigure}
\begin{subfigure}[Tensor made by random factor matrices with elements in $(-1,1)$] {\label{fig:negrand101}\includegraphics[width=0.47\textwidth, keepaspectratio]{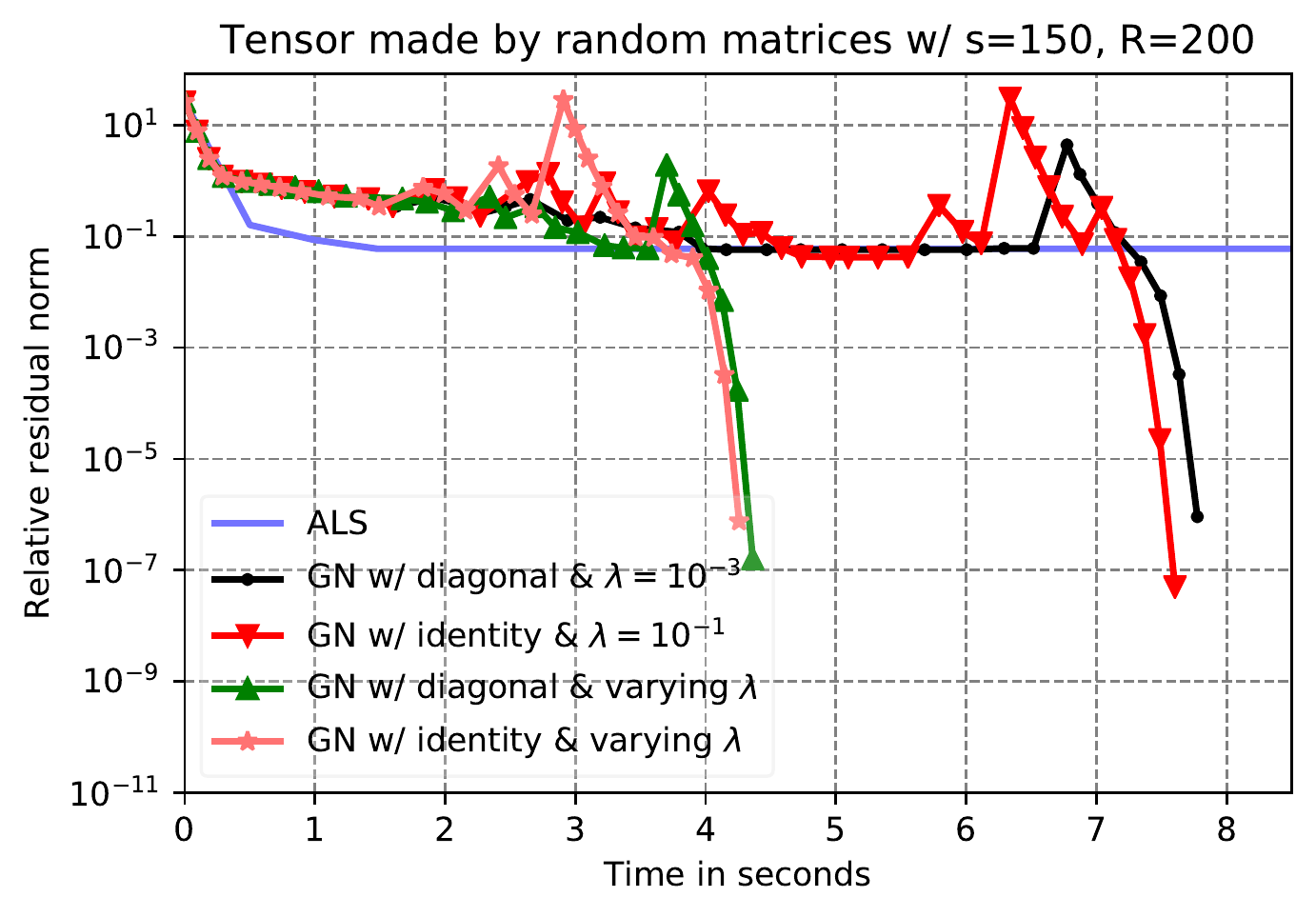}} \end{subfigure}
\caption{Relative residual norm vs time for the CP decomposition of synthetic tensors with different sizes. Timings collected using the NumPy backend (sequential).}
\label{fig:seq_behaviour}
\end{figure*}

\begin{figure*}[t]
\centering
\begin{subfigure}[Tensor made by random factor matrices with elements in $(0,1)$ using $256$ cores of Stampede2] {\label{fig:rand500}\includegraphics[width=0.48\textwidth, keepaspectratio]{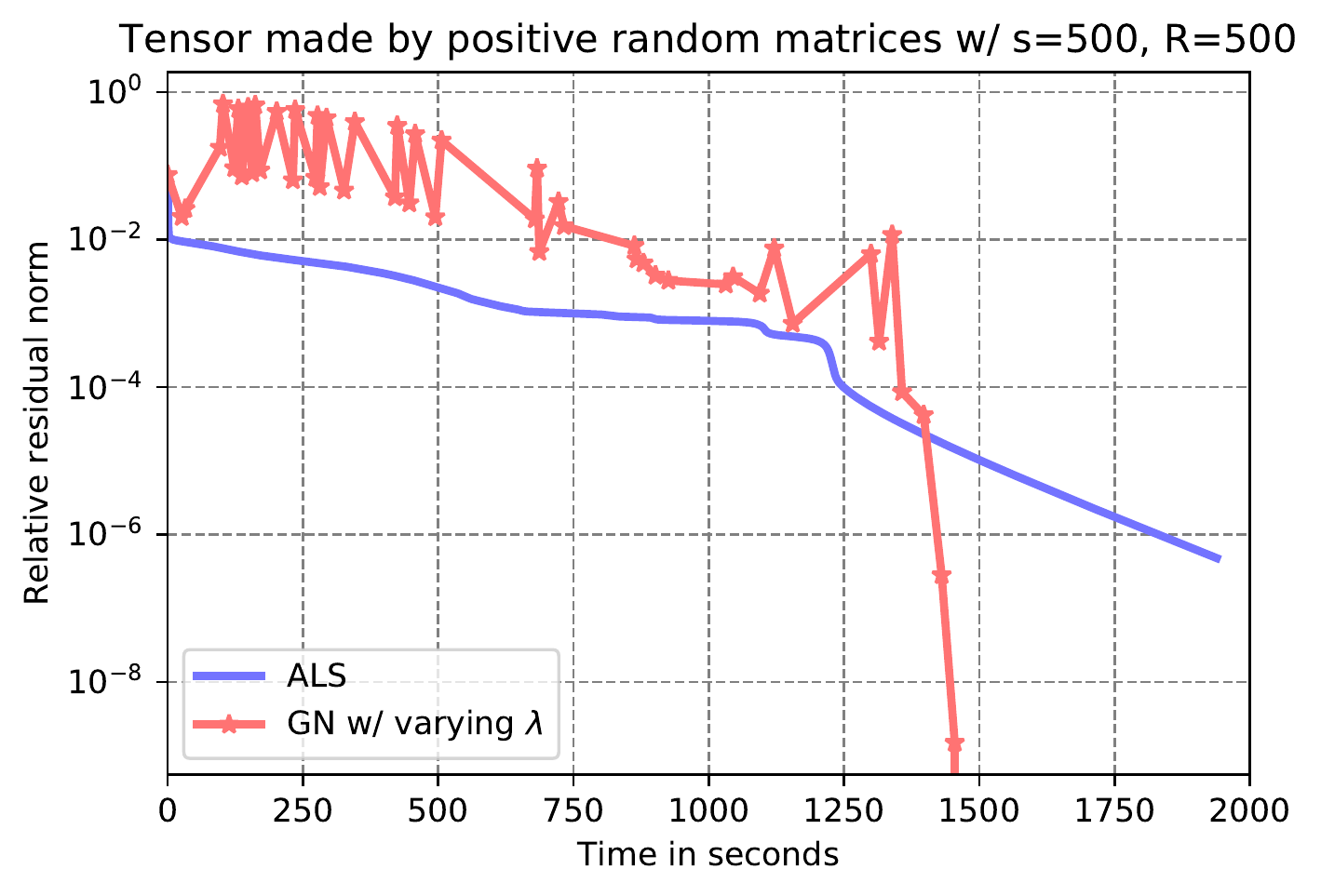}} \end{subfigure}
\begin{subfigure}[Tensor made by random factor matrices with elements distributed with standard Gaussian distribution using $256$ cores of Stampede2] {\label{fig:randn500}\includegraphics[width=0.45\textwidth, keepaspectratio]{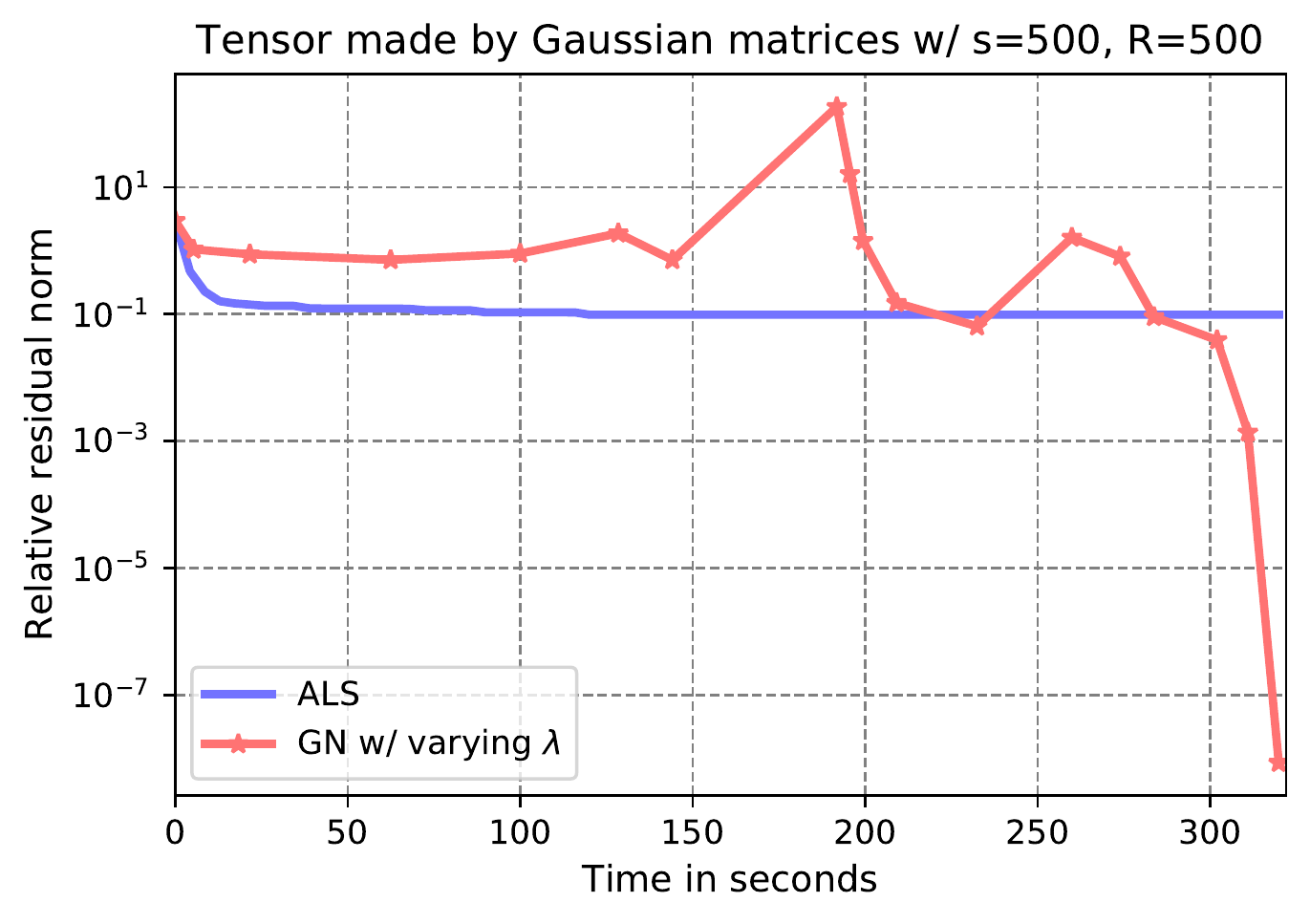}} \end{subfigure}
\begin{subfigure}[Tensor made by random factor matrices with elements in $(0,1)$ using $1024$ cores of Stampede2] {\label{fig:rand2000}\includegraphics[width=0.48\textwidth, keepaspectratio]{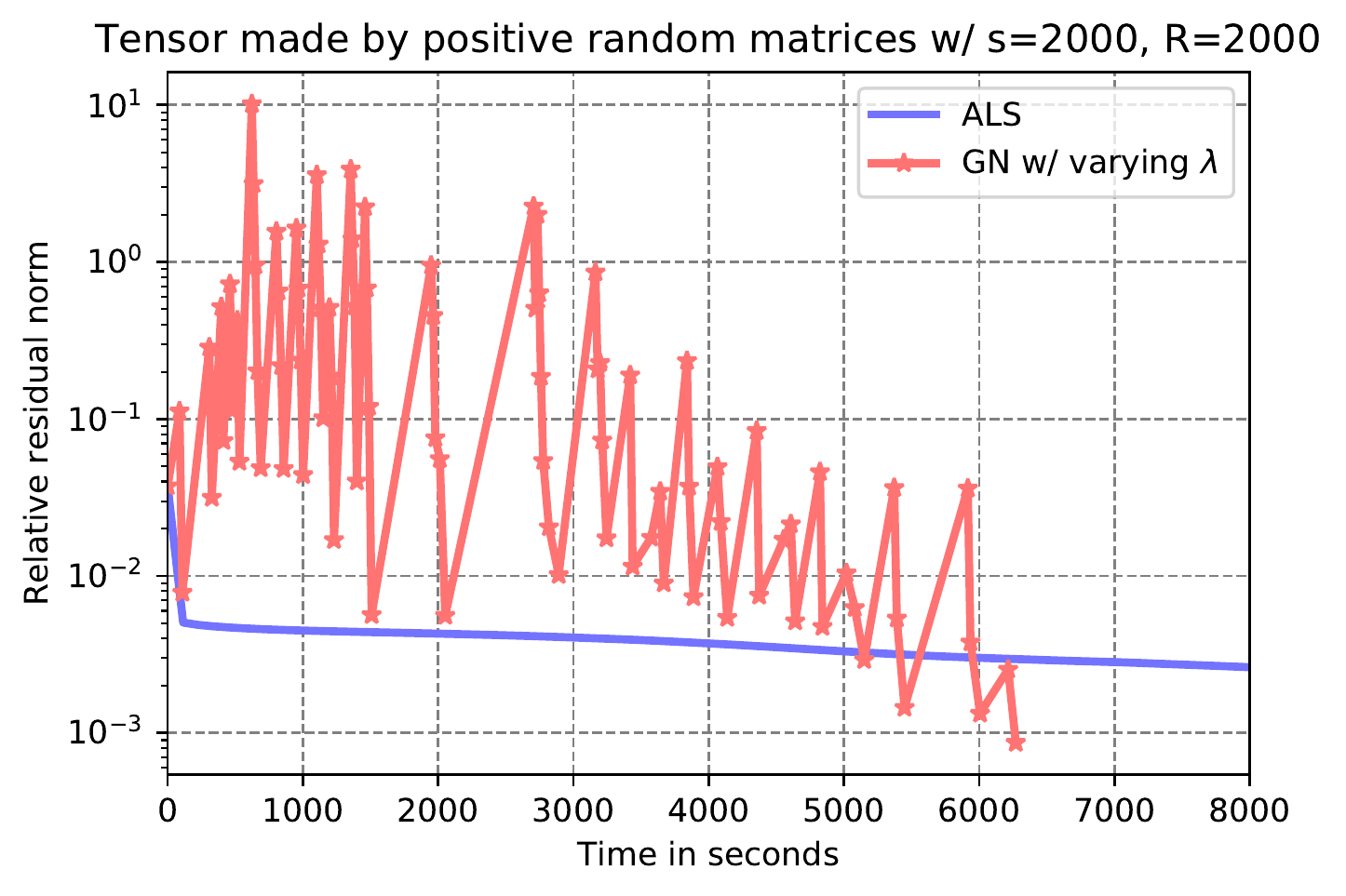}} \end{subfigure}
\caption{Relative residual norm vs time for the CP decomposition of synthetic tensors with different sizes. Timings collected using the Cyclops backend.}
\label{fig:par_behaviour}
\end{figure*}

We compare the convergence behavior of different variants of the Gauss-Newton algorithm with ALS for exact (synthetic) CP decomposition in Figure~\ref{fig:seq_behaviour} and~\ref{fig:par_behaviour}. We generate low rank tensors of different sizes, the small- and medium-sized
tensors are tested with NumPy backend and the large ones are tested with Cyclops.

In Figure~\ref{fig:gauss80} we use CP decomposition on the Gaussian low rank tensor with tensor order $N=3$, size of each dimension $s=80$ and CP rank $R=120$ with NumPy backend. We plot different types of regularization for Gauss-Newton along with ALS to study the convergence behavior of different variants of Gauss-Newton. We observe that Gauss-Newton with varying diagonal regularization performs the best and the varying identity regularization is also comparable. The sensitivity to regularization of the Gauss-Newton method is revealed in the plot as constant regularization variants are very different from each other. As can be seen in the figure, there are time periods when ALS does not make any improvement over a long time and appears to be stuck in a swamp, suggesting Gauss-Newton method is preferable for Gaussian tensors.

In Figure~\ref{fig:negrand101}, we consider the computation of CP decomposition for random low rank tensors of order three, size of each dimension $s=150$ and rank $R=200$. We can observe that constant regularization may not be useful for this tensor as we don't make any improvement over a long time. However, the other two variants with varying regularization converge fast, suggesting that varying the regularization is a robust technique for random tensors in terms of speed.

We test large random low-rank tensors in parallel with $s=500$, $R=500$ on $4$ nodes with $256$ processes as well as $s=2000$, $R=2000$ on $16$ nodes with $1024$ processes using the Cyclops backend. Gauss-Newton with identity varying regularization outperforms ALS in terms of speed and accuracy in both cases. As is shown in Figure~\ref{fig:rand500}, for $s=500$, $R=500$, Gauss-Newton with identity varying regularization converges to an exact solution about 1.25$\times$ faster than ALS, which converges to a relative residual of around $10^{-6}$. For the tensor made with standard Gaussian matrices in Figure~\ref{fig:randn500}, ALS gets stuck in a swamp
(makes very little progress in reducing the objective)
and Gauss-Newton converges to the exact solution in about 300 seconds, suggesting that for larger problems Gauss-Newton perform better for Gaussian tensors. For $s=2000$, $R=2000$ as shown in Figure~\ref{fig:rand2000}, we let the program run for a fixed time and observe Gauss-Newton with identity varying regularization converge to a lower relative residual which is about 2.4$\times$ more accurate than ALS while running for 0.6$\times$ of the time of ALS. Note that the irregularity in time taken of one Gauss-Newton iteration comes from the varying number of CG iterations taken to solve the system of equations.

\subsection{Approximate CP decomposition}
\label{subsec:approximate}

\begin{figure*}[t]
\centering

\begin{subfigure}[Input tensor size: $339\times 21\times 21$, $R=200$] {\includegraphics[width=0.48\textwidth, keepaspectratio]{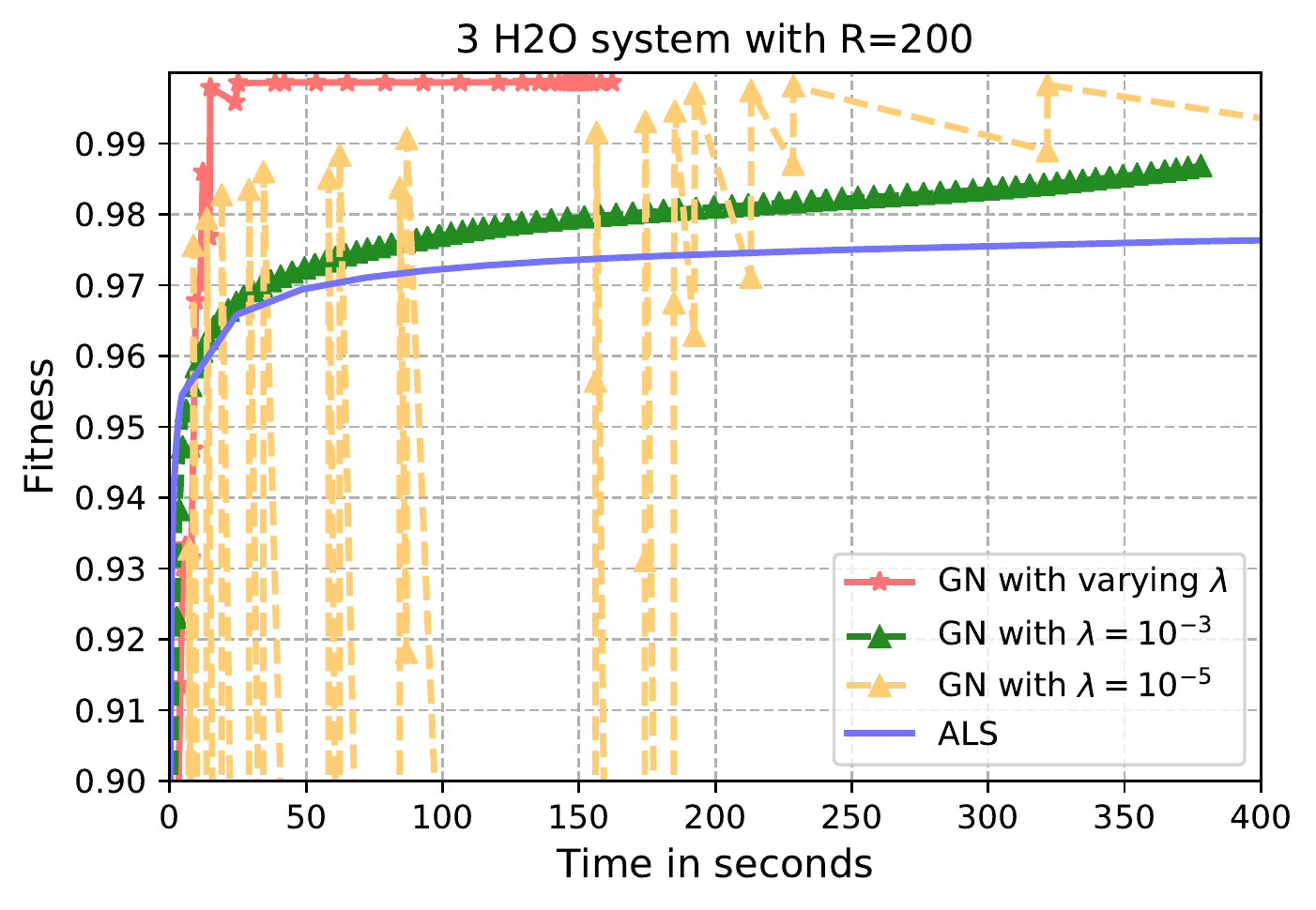}\label{subfig:quantuma}} \end{subfigure}
\begin{subfigure}[Input tensor size: $904\times 56\times 56$, $R=500$] {\label{subfig:quantumb}\includegraphics[width=0.48\textwidth, keepaspectratio]{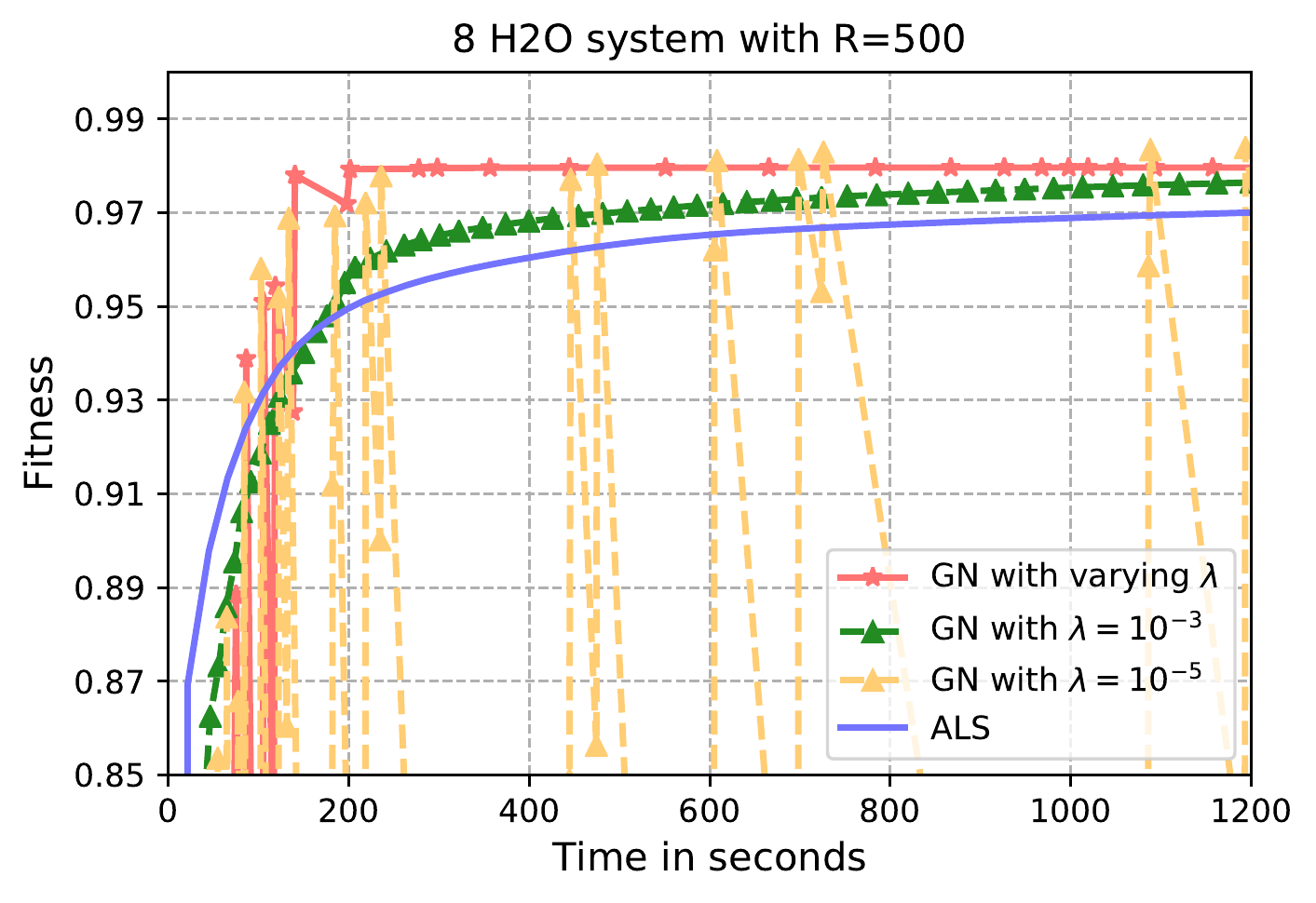}} \end{subfigure}

\begin{subfigure}[Input tensor size: $4520\times 280\times 280$, $R=2000$] {\label{subfig:quantumc}\includegraphics[width=0.48\textwidth, keepaspectratio]{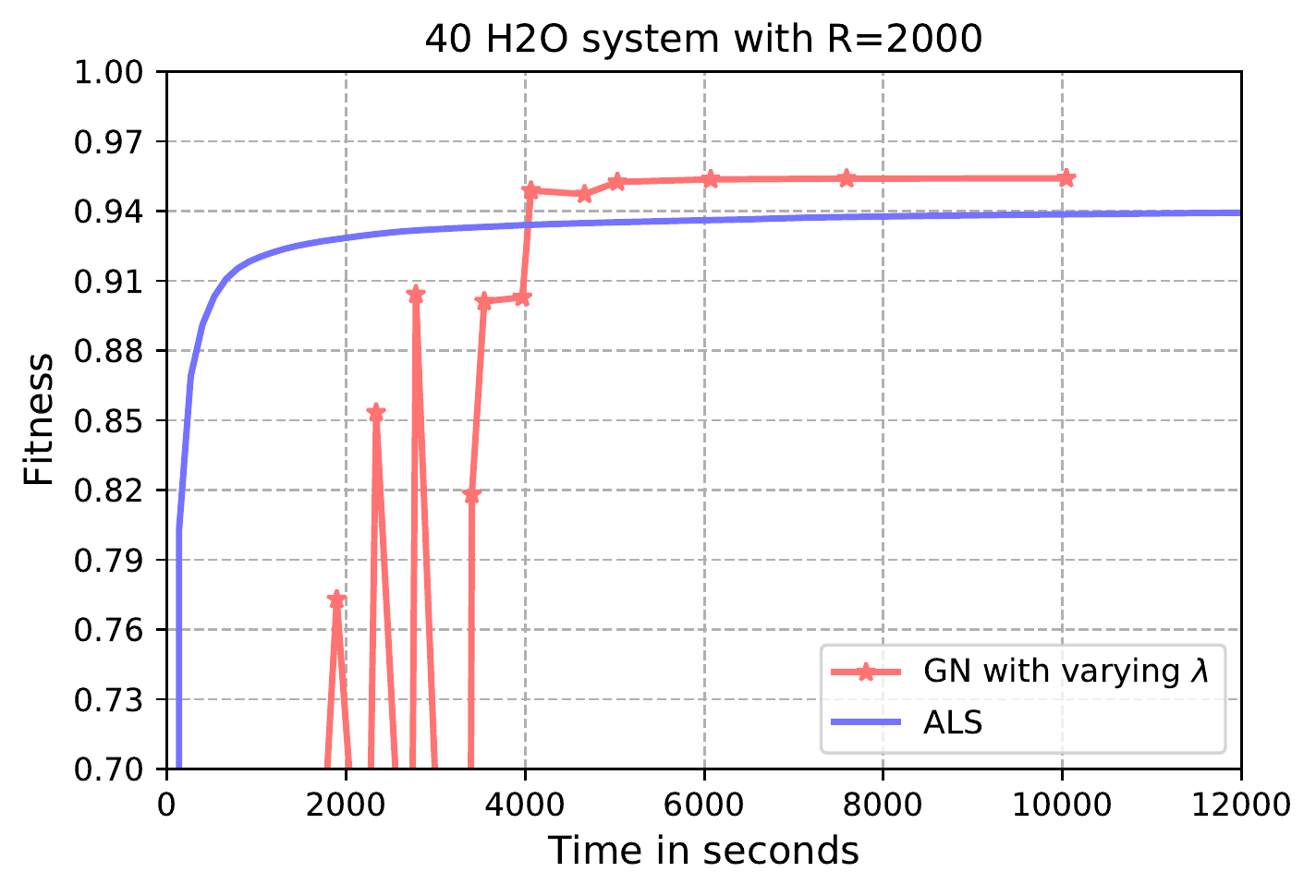}} \end{subfigure}
\begin{subfigure}[Input tensor size: $4520\times 280\times 280$, $R=3000$] {\label{subfig:quantumd}\includegraphics[width=0.48\textwidth, keepaspectratio]{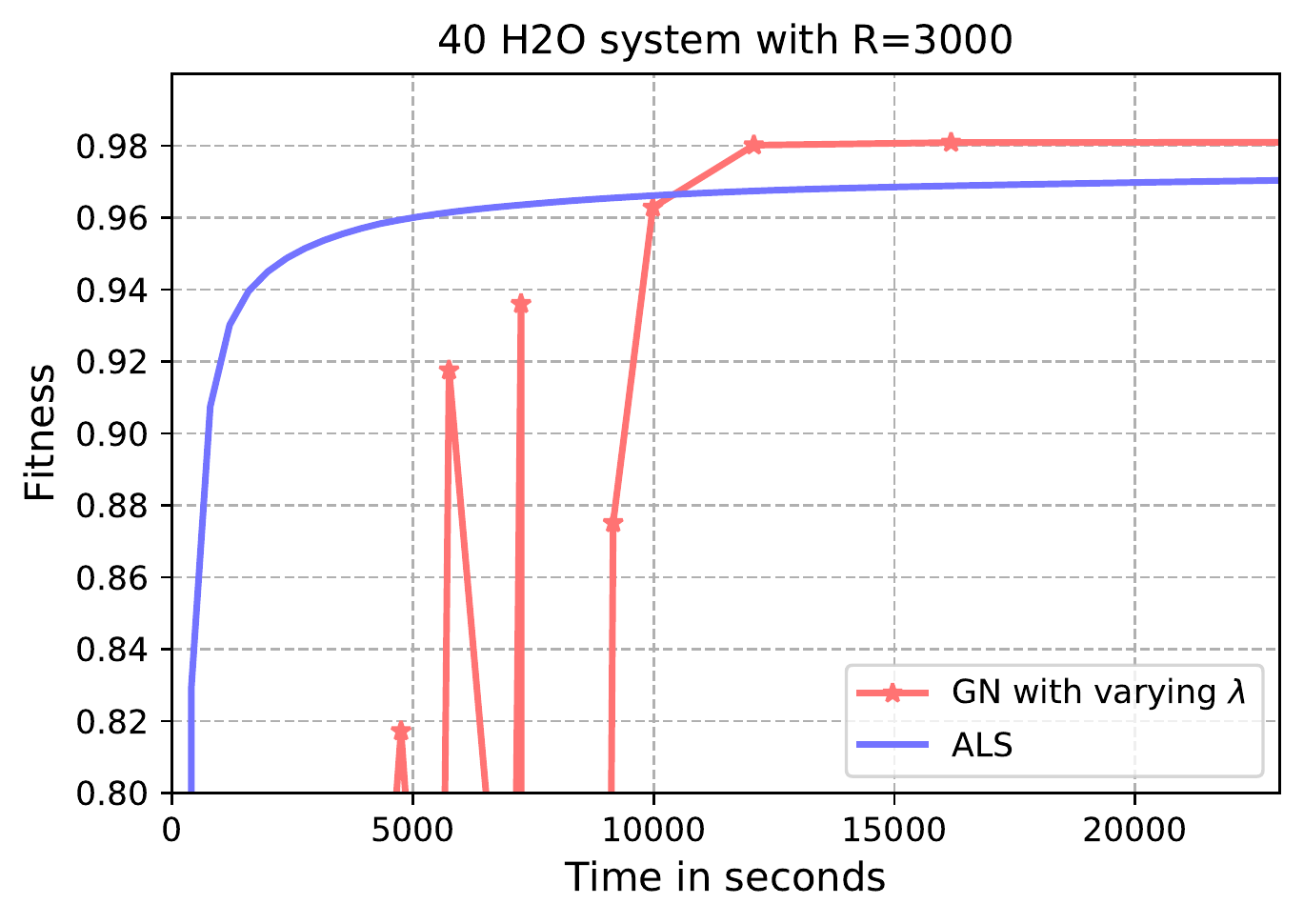}} \end{subfigure}

\caption{Fitness vs time for the CP decomposition of quantum chemistry tensors with different size and rank. The results (a) and (b) are collected with the NumPy backend, while (c) and (d) are collected with the Cyclops backend using $256$ cores of Stampede2.}
\label{fig:quantum}
\end{figure*}

We also compare the convergence behavior of Gauss-Newton method with ALS for approximate CP decomposition, in which case the tensor reconstructed from factor matrices can only approximate the input tensor rather than fully recover it. We test on the 
density fitting tensors.

Our results are shown in Figure~\ref{fig:quantum}. 
We test the problem with different input tensor sizes and different CP ranks. We run the two small sized problems shown in Figure~\ref{subfig:quantuma}, \ref{subfig:quantumb} with NumPy backend. 
We observe that for both problems, Gauss-Newton method outperforms ALS algorithm in speed and final fitness, both with the constant regularization parameter and the regularization variation scheme. In addition, Gauss-Newton with constant regularization may suffers from low optimization stability (when $\lambda=10^{-5}$) or low accuracy (when $\lambda=10^{-3}$). 
The regularization variation scheme collects the advantages of both cases, and can reach high accuracy with a stable convergence.

We run the large sized problems set up with 40 water molecules' system shown in Figure~\ref{subfig:quantumc}, \ref{subfig:quantumd} in parallel with Cyclops. Results are collected on 4 nodes using 256 processors on Stampede2. We observe that for these large problems, Gauss-Newton beats ALS in speed and fitness. With CP rank equals 2,000, Gauss-Newton can reach the fitness of 0.952 in 5,000 seconds which is higher than the best fitness of ALS (0.94) in about half the time (in 12,000 seconds), i.e., a speed up of more than 2$\times$. The oscillations in Gauss-Newton maybe controlled by using a smaller factor $\mu$ and the number of CG iterations can be reduced by using a lower regularization near the optimal solution so as to reduce the perturbation in the system of equations. The observations are similar when we increase the rank to 3000.


\section{Conclusion}
\label{sec:cnc}
In this paper, we provide the first efficient parallel implementation of a Gauss-Newton method for CP decomposition. We evaluate a formulation that employs tensor contractions for implicit matrix-vector products within the conjugate gradient method. The use of tensor contractions enables us to employ the Cyclops library for distributed-memory tensor computations to parallelize the Gauss-Newton approach with a high-level Python implementation. Our results demonstrate good weak scalability for the current implementation of the Gauss-Newton method and show how this formulation could lead to even greater speed-ups in the Hessian contraction. Additionally, we propose a regularization scheme for Gauss-Newton method to improve convergence properties without any additional cost. We perform extensive experimentation on different kinds of input tensors and compare the convergence and performance of the Gauss-Newton method relative to ALS. We observe that the Gauss-Newton method typically achieves better convergence as well as performance results for both synthetic as well as quantum chemistry tensors with high CP rank.

\section{Acknowledgments}
\label{sec:ack}

Navjot Singh, Linjian Ma, and Edgar Solomonik were supported by the US NSF OAC SSI program, award No.\ 1931258.
This work used the Extreme Science and Engineering Discovery Environment (XSEDE), which is supported by National Science Foundation grant number ACI-1548562.
We used XSEDE to employ Stampede2 at the Texas Advanced Computing Center (TACC) through allocation TG-CCR180006.
This research is part of the Blue Waters sustained-petascale computing project, which is supported by the National Science Foundation (awards OCI-0725070 and ACI-1238993) and the state of Illinois. Blue Waters is a joint effort of the University of Illinois at Urbana-Champaign and its National Center for Supercomputing Applications.

\bibliographystyle{abbrv}
\bibliography{paper}

\end{document}